\documentclass[english]{article}
\usepackage[T1]{fontenc}
\usepackage[latin9]{inputenc}
\usepackage{verbatim}
\usepackage{amsmath}
\usepackage{amssymb}
\usepackage{geometry}
\usepackage{esint}

\makeatletter
\newcommand{\lyxaddress}[1]{
	\par {\raggedright #1
	\vspace{1.4em}
	\noindent\par}
}

\makeatother

\makeatother

\usepackage{babel}
\begin{document}
\title{An infinite set of one-range addition theorems without an infinite
second series, for Slater orbitals and it derivatives, applicable more
than one coordinate system}
\author{Jack C. Straton}
\maketitle

\lyxaddress{Department of Physics, Portland State University, Portland, OR, 97207-0751,
straton@pdx.edu}
\begin{abstract}
$^{\text{}}$Addition theorems have been indispensable tools for the
reduction of quantum transition amplitudes. They are normally utilized
at the start of the process to move the angular dependence within
plane waves and Coulomb potentials, and the like, into a sum over
Spherical Harmonics that allows the angular integration to be carried
out. These have historically been ``two-range'' addition theorems,
characterized by the two-fold notation $r_{>}=Max[r_{1},r_{2}]$ and
$r_{<}=Min[r_{1},r_{2}]$ and comprising a single infinite series.
More recently, ``one-range'' addition theorems have been created
that have no such piecewise notation, but at the cost of a second
infinite series. We use a very different approach to derive an infinite
set of addition theorems for Slater orbitals and its derivatives that
retain the one-range variable dependence but have, at worst, a finite
second series rather than an infinite one. Also unlike previous addition
theorems, they are applicable to more than one coordinate system.
One of these addition theorem may also be used for Yukawa-like functions
that may appear late in the reduction of amplitude integrals and we
show its utility for an integral that has stubbornly defied reduction
to analytic form for nearly sixty years. 
\end{abstract}

\vspace{2pc}
 \\
MSC classes: 44A20, 44A30,  81Q99, 33C10, 33C20,  33C60, 34B27

\section{Introduction}

The central impediment to reducing the dimensionality of quantum transition
amplitudes is the presence of angular cross-terms $|{\bf x}_{1}-{\bf x}_{2}|=\sqrt{x_{1}^{2}-2x_{1}x_{2}\cos\theta+x_{2}^{2}}$
sequestered in various square roots of quadratic forms. Direct spatial
integration is sometimes possible (see, for instance, \cite{Ley-Koo and Bunge},
among many others), and at other times addition theorems (e.g., \cite{Sack},
\cite{Porras and King}, \cite{Weniger_two-range}, and \cite{Guseinov2005})
are more useful. One may instead apply Fourier transforms (e.g., \cite{Fromm and Hill},
\cite{Remiddi}, and \cite{Harris PRA 55 1820}) and/or Gaussian transforms
(e.g., \cite{Kikuchi}, \cite{Shavitt and Karplus}, and \cite{Stra89a})
to effect these reductions.

The author has recently introduced \cite{stra23} a fifth reduction
method in the spirit of Fourier and Gaussian transforms that constitutes
an M-1-dimensional integral representation over the interval $[0,\infty]$
for products of M Slater orbitals, which is one fewer integral dimension
than a Gaussian transform requires, and roughly a quarter of the integral
dimensions introduced by a Fourier transform for such a product. Subsequent
work constructed \cite{stra24a} a sixth reduction method using an
M-1-dimensional integral representation over the interval $[0,1]$.

A comparison of all of these approaches often centers on the Slater
orbital $\frac{e^{-\eta x}}{x}$, commonly written without arguments
as $\psi_{000}$, which acts as a seed function from which Slater
functions,\cite{Chen} Hylleraas powers,\cite{Harris PRA 55 1820}
and hydrogenic wave functions can be derived by differentiation. In
nuclear physics it is known as the Yukawa\cite{Yukawa} exchange potential,
and in plasma physics as the Debye-H\"uckel  potential, arising from
screened charges \cite{NayekGhoshal} so that the Coulomb potential
is replaced by an effective screened potential.\cite{Harris, EckerWeizel}
Screening of charges also appears in solid-state physics, where this
function is called the Thomas-Fermi potential. In the atomic physics
of negative ions, the radial wave function is given by the equivalent
Macdonald function $\left(R(r)=\frac{C}{\sqrt{r}}K_{1/2}(\eta r)\right)$.
\cite{Smirnov2003} This function also appears in the approximate
ground state wave function \cite{GaravelliOliveira} for a hydrogen
atom interacting with hypothesized non-zero-mass photons.\cite{CaccavanoLeung}
With imaginary $\eta$, this function is also ($4\pi$ times) the
free-space Green\textquoteright s function, (see, for instance, \cite{Stilgoe},
among many others). We will simply call these \emph{Slater orbitals}
herein.

The present work focuses on addition theorems. Probably the most familiar
addition theorem, for the Coulomb potential \cite{Joachain} (p. 670,
Eq. (B.40)), is

\begin{equation}
\psi_{000}\left(0,x_{12}\right)=\frac{1}{x_{12}}\equiv\frac{1}{\left|\mathbf{r}_{1}-\mathbf{r}_{2}\right|}=4\pi\sum_{l=0}^{\infty}\sum_{m=-l}^{l}\frac{1}{2l+1}\frac{r_{<}^{l}}{r_{>}^{l+1}}Y_{lm}^{*}(\theta_{1},\phi_{1})Y_{lm}(\theta_{2},\phi_{2})\:.\label{eq:Coulomb}
\end{equation}
 This is an example of a ``two-range addition theorem'' characterized
by the two-fold notation $r_{>}=Max[r_{1},r_{2}]$ and $r_{<}=Min[r_{1},r_{2}]$
and comprising a single infinite series. It may be utilized in the
simplest example of a transition amplitude so that \cite{Joachain}
(p. 669 (B.35))

\begin{eqnarray}
S_{1}^{\eta_{1}000}\left(0;0,\mathbf{x}_{2}\right) & \equiv & S_{1}^{\eta_{1}j_{1}\eta_{2}j_{2}}\left(\mathbf{p}_{1};\mathbf{y}_{1},\mathbf{y}_{2}\right)_{p_{1}\rightarrow0,y_{1},\rightarrow0,y_{2}\rightarrow x_{2},j_{1}\rightarrow0,\eta_{2}\rightarrow0,j_{2}\rightarrow0}\nonumber \\
 & = & \int d^{3}x_{1}\frac{e^{-\eta_{1}x_{1}}}{x_{1}}\frac{1}{x_{12}}\nonumber \\
 & = & \int_{0}^{\infty}x_{1}^{2}dx_{1}\frac{e^{-\eta_{1}x_{1}}}{x_{1}}\int d\Omega_{1}\left[\sqrt{4\pi}Y_{00}\left(\Omega_{1}\right)\right]4\pi\sum_{l=0}^{\infty}\sum_{m=-l}^{l}\frac{1}{2l+1}\frac{r_{<}^{l}}{r_{>}^{l+1}}Y_{lm}^{*}\left(\Omega_{1}\right))Y_{lm}\left(\Omega_{2}\right)\nonumber \\
 & = & 4\pi\int_{0}^{\infty}x_{1}^{2}dx_{1}\frac{e^{-\eta_{1}x_{1}}}{x_{1}}\sum_{l=0}^{\infty}\sum_{m=-l}^{l}\frac{1}{2l+1}\frac{r_{<}^{l}}{r_{>}^{l+1}}\delta_{l0}\delta_{lm0}\sqrt{4\pi})Y_{lm}\left(\Omega_{2}\right)\nonumber \\
 & = & 4\pi\left(\int_{0}^{x_{2}}x_{1}^{2}dx_{1}\frac{e^{-\eta_{1}x_{1}}}{x_{1}}\frac{1}{x_{2}}+\int_{x_{2}}^{\infty}x_{1}^{2}dx_{1}\frac{e^{-\eta_{1}x_{1}}}{x_{1}}\frac{1}{x_{1}}\right)\nonumber \\
 & = & \frac{4\pi\left(1-e^{-\eta_{1}x_{2}}\right)}{x_{2}\eta_{1}^{2}}\:,\label{eq:syc}\\
 &  & \qquad\left[x_{2}>0,\,Re\,\eta_{1}>0\right]\nonumber 
\end{eqnarray}
 where the notation $r_{>}$ is made manifest in the second to last
line. Here we use the much more general notation of previous work
\cite{Stra89a} in which the short-hand form for shifted coordinates
is $\mathbf{x}_{12}=\mathbf{x}_{1}-\mathbf{x}_{2}$, $\mathbf{p}_{1}$
is a momentum variable within any plane wave associated with the (first)
integration variable, the $\mathbf{y}_{i}$ are coordinates external
to the integration, and the \emph{j}s are defined in the Gaussian
transform \cite{Stra89a} of the generalized Slater orbital (\ref{seventeen}),
below. One can also take derivatives of the Slater orbital $e^{-\eta_{1}x_{1}}/x_{1}$
with respect to $\eta_{1}$ in the integrand to obtain functions with
$j>0$, such as hydrogenic s-states, and likewise of the solution
on the last line.

Starting with functions having the Slater orbital form gives fewer
terms at each integral-reduction step than does starting with hydrogenic
wave functions.

Another often-used addition theorem is the one for plane waves \cite{Joachain}
(p. 671 eq. (B.44))
\begin{equation}
e^{i\mathbf{k\cdot r}}=4\pi{\displaystyle \sum_{L=0}^{\infty}}{\displaystyle \sum_{M=-L}^{L}}i^{L}j_{L}(kr)Y_{LM}^{*}\left(\theta_{k},\phi_{k}\right)Y_{LM}\left(\theta,\phi\right)\label{eq:planewave}
\end{equation}
 that is useful in problems like photoionization \cite{Casey} in
which the interaction potential takes on the form $r\cos\theta$.

If both of the functions have the Slater orbital form, one may use
a lesser-known addition theorem found in Magnus, Oberhettinger, and
Soni \cite{Magnus Oberhettinger and Soni}

\begin{equation}
\begin{array}{ccc}
{\displaystyle \psi_{000}\left(\eta,x_{12}\right)=\frac{e^{-\eta\sqrt{x_{1}^{2}-2\mathbf{x}_{1}\cdot\mathbf{x}_{2}+x_{2}^{2}}}}{\sqrt{x_{1}^{2}-2\mathbf{x}_{1}\cdot\mathbf{x}_{2}+x_{2}^{2}}}} & = & x_{1}^{-1/2}x_{2}^{-1/2}{\displaystyle \sum_{n=0}^{\infty}}\left(2n+1\right)P_{n}(\cos(\theta))I_{n+1/2}\left(\eta x_{1}\right)K_{n+1/2}\left(\eta x_{2}\right)\\
 &  & \left[0<x_{1}<x_{2}\right]
\end{array}\label{eq:yukseriesx1x2-1}
\end{equation}
 that may be cast into modern notation \cite{Joachain} (p. 670 (B.42)
has the equivalent version with imaginary $\eta$)

\begin{equation}
\psi_{000}\left(\eta,x_{12}\right)=\begin{array}{ccc}
{\displaystyle \frac{e^{-\eta\sqrt{x_{1}^{2}-2\mathbf{x}_{1}\cdot\mathbf{x}_{2}+x_{2}^{2}}}}{\sqrt{x_{1}^{2}-2\mathbf{x}_{1}\cdot\mathbf{x}_{2}+x_{2}^{2}}}} & = & x_{1}^{-1/2}x_{2}^{-1/2}{\displaystyle \sum_{n=0}^{\infty}}\left(2n+1\right)P_{n}(\cos(\theta))I_{n+1/2}\left(\eta x_{<}\right)K_{n+1/2}\left(\eta x_{>}\right)\end{array} \;.\label{eq:yukseriesx1x2><}
\end{equation}

\subsection{The central problem with two-range addition theorems }

This allows one to solve the next most complicated integral, \cite{GR5}
( p. 155 No. 2.481.1 with $\eta_{1}^{2}\neq\eta_{2}^{2}$, and GR5
p. 357 No. 3.351.2)

\begin{eqnarray}
S_{1}^{\eta_{1}0\eta_{2}0}\left(0;0,\mathbf{x}_{2}\right) & = & \int d^{3}x_{1}\frac{e^{-\eta_{1}x_{1}}}{x_{1}}\frac{e^{-\eta_{2}x_{12}}}{x_{12}}\nonumber \\
 & = & \int_{0}^{\infty}dx_{1}x_{1}^{2-1-1/2}e^{-\eta_{1}x_{1}}x_{2}^{-1/2}{\displaystyle \sum_{n=0}^{\infty}}\left(2n+1\right)I_{n+1/2}\left(\eta_{2}x_{<}\right)K_{n+1/2}\left(\eta_{2}x_{>}\right)\nonumber \\
 & \times & 2\pi\int_{-\pi}^{\pi}d(\cos(\theta))P_{0}(\cos(\theta))P_{n}(\cos(\theta))\nonumber \\
 & = & x_{2}^{-1/2}\int_{0}^{\infty}dx_{1}x_{1}^{2-1-1/2}e^{-\eta_{1}x_{1}}{\displaystyle \sum_{n=0}^{\infty}}\left(2n+1\right)I_{n+1/2}\left(\eta_{2}x_{<}\right)K_{n+1/2}\left(\eta_{2}x_{>}\right)\nonumber \\
 & \times & 2\pi\int_{-1}^{1}duP_{0}(u)P_{n}(u)\nonumber \\
 & = & x_{2}^{-1/2}\int_{0}^{\infty}dx_{1}x_{1}^{2-1-1/2}e^{-\eta_{1}x_{1}}{\displaystyle \sum_{n=0}^{\infty}}\left(2n+1\right)I_{n+1/2}\left(\eta_{2}x_{<}\right)K_{n+1/2}\left(\eta_{2}x_{>}\right)\nonumber \\
 & \times & 2\pi\frac{2}{2\cdot0+1}\delta_{0n}\nonumber \\
 & = & 4\pi x_{2}^{-1/2}\left(K_{1/2}\left(\eta_{2}x_{2}\right)\int_{0}^{x_{2}}dx_{1}x_{1}^{1/2}e^{-\eta_{1}x_{1}}I_{1/2}\left(\eta_{2}x_{1}\right)\right.\nonumber \\
 & + & \left.\hspace{1.6cm}I_{1/2}\left(\eta_{2}x_{2}\right)\int_{x_{2}}^{\infty}dx_{1}x_{1}^{1/2}e^{-\eta_{1}x_{1}}K_{1/2}\left(\eta_{2}x_{1}\right)\right)\nonumber \\
 & = & 4\pi x_{2}^{-1/2}\left(\sqrt{\frac{\pi}{2\eta_{2}x_{2}}}e^{-\eta_{2}x_{2}}\int_{0}^{x_{2}}dx_{1}x_{1}^{1/2}e^{-\eta_{1}x_{1}}\frac{\sqrt{\frac{2}{\pi}}\sinh\left(x_{1}\eta_{2}\right)}{\sqrt{x_{1}\eta_{2}}}\right.\nonumber \\
 & + & \left.\hspace{1.7cm}\frac{\sqrt{\frac{2}{\pi}}\sinh\left(x_{2}\eta_{2}\right)}{\sqrt{x_{2}\eta_{2}}}\int_{x_{2}}^{\infty}dx_{1}x_{1}^{1/2}e^{-\eta_{1}x_{1}}\sqrt{\frac{\pi}{2\eta_{2}x_{1}}}e^{-\eta_{2}x_{1}}\right)\nonumber \\
 & = & \frac{4\pi}{x_{2}}\left(\frac{e^{-x_{2}\eta_{1}-x_{2}\eta_{2}}\left(\eta_{1}\sinh\left(x_{2}\eta_{2}\right)+\eta_{2}\left(\cosh\left(x_{2}\eta_{2}\right)-e^{x_{2}\eta_{1}}\right)\right)}{\eta_{2}\left(\eta_{2}^{2}-\eta_{1}^{2}\right)}\right.\nonumber \\
 & + & \left.\hspace{3.8cm}\frac{\left(\eta_{2}-\eta_{1}\right)e^{-x_{2}\left(\eta_{1}+\eta_{2}\right)}\sinh\left(x_{2}\eta_{2}\right)}{\eta_{2}\left(\eta_{2}^{2}-\eta_{1}^{2}\right)}\right)\nonumber \\
 & = & \frac{4\pi\left(e^{-\eta_{2}x_{2}}-e^{-\eta_{1}x_{2}}\right)}{x_{2}\left(\eta_{1}^{2}-\eta_{2}^{2}\right)}\;,\label{eq:SVxVxyG}
\end{eqnarray}
 which is correct.

While this approach seems to be an excellent path forward, for more
complicated integrals it manifests problems. For instance, the author
attempted to use this approach for a product of five Slater orbitals
with shifted coordinates, and split integrals in this two-range approach
manifested pairs of large-amplitude, nearly-canceling, terms in numerical
checks of each step that eventually brought the project to a standstill
for lack of accuracy.

One-range addition theorems (see for instance \cite{Guseinov2005})
remove the need for such notations. More crucially, they also bypass
significant negative numerical consequences of round-off errors in
cancelling terms of large magnitude and opposite signs that may appear
in applications of two-range addition theorems. But one-range addition
theorems do so at the cost of an additional infinite series:

\begin{equation}
\begin{array}{ccc}
\psi_{000}\left(\zeta,x_{12}\right)={\displaystyle \frac{e^{-\zeta\sqrt{x_{1}^{2}-2\mathbf{x}_{1}\cdot\mathbf{x}_{2}+x_{2}^{2}}}}{\sqrt{x_{1}^{2}-2\mathbf{x}_{1}\cdot\mathbf{x}_{2}+x_{2}^{2}}}} & = & \hspace{-4.8cm}\frac{2^{3/2}}{\left(2\zeta\right)^{2}}\underset{\begin{array}{c}
N\rightarrow\infty\\
N'\rightarrow\infty
\end{array}}{\lim}{\displaystyle \sum_{n=0}^{N}}{\displaystyle \sum_{l=0}^{n-1}}{\displaystyle \sum_{m=-l}^{l}}\\
 &  & \left[{\displaystyle \sum_{u=0}^{N'}}{\displaystyle \sum_{v=0}^{u-1}}{\displaystyle \sum_{s=-v}^{v}}B_{000,nlm}^{\alpha uvs}\left(N,N';\eta,\zeta\right)\chi_{uvs}^{*}\left(\zeta,\mathbf{x}_{2}\right)\right]\chi_{nlm}\left(\zeta,\mathbf{x}_{1}\right)
\end{array}\label{eq:one-range}
\end{equation}
 where $\alpha=1,0,-1,-2,\cdots$, the Slater-Type Orbitals (STOs)
are

\begin{equation}
\chi_{\mu\nu\sigma}\left(\zeta,\mathbf{x}_{1}\right)=\frac{\left(2\zeta\right)^{\mu+1/2}}{\left(2\mu\right)!}r^{\mu-1}e^{-\zeta x}Y_{\nu\sigma}\left(\theta,\phi\right)\:,\label{eq:STO}
\end{equation}
 and B is a complicated function of its parameters given in \cite{Guseinov2005})
and \cite{Guseinov2002}.

\section{Integral Transforms}

Integral transforms provide the main alternative to addition theorems
for transition amplitudes. Consider, for instance, the Fourier transform
of a product of Slater orbitals, 
\begin{equation}
S_{1}^{\eta_{1}0\eta_{2}0}\left(k,;0,x_{2}\right)=\int d^{3}x_{1}\frac{e^{-\eta_{1}x_{1}}}{x_{1}}\frac{e^{-\eta_{2}x_{12}}}{x_{12}}e^{-i\mathbf{k}\cdot\mathbf{x}_{1}}\label{eq:kyy}
\end{equation}
 The Slater orbital and its derivatives have Gaussian transform \cite{Stra89a}
\begin{equation}
\begin{array}[t]{ccc}
V^{\eta j}({\bf R}) & = &R^{j-1}e^{-\eta R}=\left(-1\right)^{j}{\displaystyle \frac{d^{j}}{d\eta^{j}}}  {\displaystyle \frac{1}{\sqrt{\pi}}}\int_{0}^{\infty}\,d\rho_{3}{\displaystyle \frac{e^{-R^{2}\rho_{3}}e^{-\eta^{2}/4/\rho_{3}}}{\rho_{3}^{\;1/2}}}\;\;\left[\eta\geqq0,\:R>0\right]\\
 & =R^{j-1}e^{-\eta R}= & {\displaystyle \frac{1}{2^{j}\sqrt{\pi}}}\int_{0}^{\infty}\,d\rho_{3}{\displaystyle \frac{e^{-R^{2}\rho_{3}}e^{-\eta^{2}/4/\rho_{3}}}{\rho_{3}^{\;\left(j+1\right)/2}}}H_{j}\left({\displaystyle \frac{\eta}{2\sqrt{\rho_{3}}}}\right)\;\;\left[\forall j\:\mathrm{if\:}\eta>0,\;j=2J\:\mathrm{if\:}\eta=0\right]\\
 &  & \left[\forall j\:\mathrm{if\:}\eta>0,\;j=2J\:\mathrm{if\:}\eta=0,\textrm{ integrate first if }j=2J+1\:\mathrm{and\:}\eta=0\right]
\end{array}\label{seventeen}
\end{equation}
 that allow us to combine all angular dependence into one quadratic
form, and then complete the square.

\begin{eqnarray}
S_{1}^{\eta_{1}0\eta_{2}0}\left(\mathbf{k};0,\mathbf{x}_{2}\right)\  & = & \int d^{3}x_{1}{\displaystyle \frac{1}{\sqrt{\pi}}}\int_{0}^{\infty}\,d\rho_{1}{\displaystyle \frac{e^{-x_{1}^{2}\rho_{1}}e^{-\eta_{1}^{2}/4/\rho_{1}}}{\rho_{1}^{\;1/2}}}{\displaystyle \frac{1}{\sqrt{\pi}}}\int_{0}^{\infty}\,d\rho_{2}{\displaystyle \frac{e^{-x_{12}^{2}\rho_{2}}e^{-\eta_{2}^{2}/4/\rho_{2}}}{\rho_{2}^{\;1/2}}e^{-i\mathbf{k}\cdot\mathbf{x}_{1}}}\nonumber \\
 & = & {\displaystyle \frac{1}{\pi}}\int d^{3}x'_{1}\int_{0}^{\infty}\,d\rho_{1}{\displaystyle \frac{e^{-\eta_{1}^{2}/4/\rho_{1}}}{\rho_{1}^{\;1/2}}}\int_{0}^{\infty}\,d\rho_{2}{\displaystyle \frac{e^{-\eta_{2}^{2}/4/\rho_{2}}}{\rho_{2}^{\;1/2}}}\nonumber \\
 & \times & exp\left(-\left(\rho_{1}+\rho_{2}\right)x'{}_{1}^{2}-\frac{x_{2}^{2}\rho_{1}\rho_{2}}{\rho_{1}+\rho_{2}}-\frac{i\mathbf{k}\cdot\mathbf{x}_{2}\rho_{2}}{\rho_{1}+\rho_{2}}-\frac{\text{k}^{2}}{4\left(\rho_{1}+\rho_{2}\right)}\right)\:,\label{eq:squarecompleted}
\end{eqnarray}
 where we have changed variables from $\mathbf{x}_{1}$ to $\mathbf{x}'_{1}=\mathbf{x}_{1}-\frac{\mathbf{x}_{2}\rho_{2}-i\mathbf{k}/2}{\rho_{1}+\rho_{2}}$
with unit Jacobian so that the spatial integral may be done: \cite{GR5}
(p. 382 No. 3.461.2)

\begin{equation}
\int e^{-\left(\rho_{1}+\rho_{2}\right)x{'}_{1}^{2}}d^{3}x{'}_{1}=4\pi\int_{0}^{\infty}e^{-\left(\rho_{1}+\rho_{2}\right)x{'}_{1}^{2}}\;x{'}_{1}^{2}dx'_{1}=\frac{4\pi^{1+1/2}}{2^{2}\left(\rho_{1}+\rho_{2}\right)^{3/2}}\quad\left[\rho_{1}+\rho_{2}>0\right].\label{eq:spatial integral}
\end{equation}

What is left over is

\begin{eqnarray}
S_{1}^{\eta_{1}0\eta_{2}0}\left(\mathbf{k};0,\mathbf{x}_{2}\right) & = & \pi^{1/2}\int_{0}^{\infty}\,d\rho_{1}{\displaystyle \frac{e^{-\eta_{1}^{2}/4/\rho_{1}}}{\rho_{1}^{\;1/2}}}\int_{0}^{\infty}\,d\rho_{2}{\displaystyle \frac{e^{-\eta_{2}^{2}/4/\rho_{2}}}{\rho_{2}^{\;1/2}}}\label{eq:rho_version}\\
 & \times & \frac{1}{\left(\rho_{1}+\rho_{2}\right)^{3/2}}\exp\left(-\frac{x_{2}^{2}\rho_{1}\rho_{2}}{\rho_{1}+\rho_{2}}-\frac{i\mathbf{k}\cdot\mathbf{x}_{2}\rho_{2}}{\rho_{1}+\rho_{2}}-\frac{\text{k}^{2}}{4\left(\rho_{1}+\rho_{2}\right)}\right)\:.\nonumber 
\end{eqnarray}

\noindent Let

%\begin{equation}
%\begin{array}{cc}
%\tau_{1}={\displaystyle \frac{\rho_{2}}{\rho_{1}+\rho_{2}}},\qquad\rho_{1}=\rho_{2}\,{\displaystyle \frac{\left(1-\tau_{1}\right)}{\tau_{1}}},\rho_{2}=\rho_{1}{\displaystyle \frac{\tau_{1}}{\left(1-\tau_{1}\right)}},{\displaystyle \frac{1}{\rho_{1}+\rho_{2}}},=\frac{\tau_{1}}{\rho_{2}}=\frac{\left(1-\tau_{1}\right)}{\rho_{1}}\\
%d\tau_{1}=d\rho_{2}\left({\displaystyle \frac{1}{\rho_{1}+\rho_{2}}-{\displaystyle \frac{\rho_{2}}{\left[\rho_{1}+\rho_{2}\right]^{2}}}}\right)=d\rho_{2}\left({\displaystyle {\displaystyle \frac{\rho_{2}}{\left[\rho_{1}+\rho_{2}\right]^{2}}}}\right)\\
%d\tau_{1}=d\rho_{2}\left({\displaystyle {\displaystyle \frac{\rho_{1}\tau_{1}^{2}}{\rho_{2}^{2}}}}\right)=d\rho_{2}\left({\displaystyle {\displaystyle {\displaystyle {\displaystyle \frac{\tau_{1}^{2}}{\rho_{2}^{2}}}}}\rho_{2}\,{\displaystyle \frac{\left(1-\tau_{1}\right)}{\tau_{1}}}}\right)=d\rho_{2}\left({\displaystyle {\displaystyle \frac{\tau_{1}\left(1-\tau_{1}\right)}{\rho_{2}}}}\right)\:.\\
%d\tau_{1}={\displaystyle \frac{d\rho_{2}}{\rho_{2}^{\;1/2}}}\left({\displaystyle {\displaystyle \frac{\tau_{1}\left(1-\tau_{1}\right)}{\rho_{1}^{1/2}{\displaystyle \frac{\tau_{1}^{1/2}}{\left(1-\tau_{1}\right)^{1/2}}}}}}\right)={\displaystyle \frac{d\rho_{2}}{\rho_{2}^{\;1/2}}}\left({\displaystyle {\displaystyle \frac{\tau_{1}^{1/2}\left(1-\tau_{1}\right)^{3/2}}{\rho_{1}^{1/2}}}}\right)\end{array}\label{eq:tau transform}\end{equation}

\begin{equation}
\tau_{1}={\displaystyle \frac{\rho_{2}}{\rho_{1}+\rho_{2}}},\quad d\tau_{1}={\displaystyle \frac{d\rho_{2}}{\rho_{2}^{\;1/2}}}\left({\displaystyle {\displaystyle \frac{\tau_{1}^{1/2}\left(1-\tau_{1}\right)^{3/2}}{\rho_{1}^{1/2}}}}\right)\ \label{eq:tau transform}
\end{equation}

\noindent When $\rho_{2}=0$ then $\tau_{1}=0$ and when $\rho_{2}=\infty$
then $\tau_{1}=1$ . Then \cite{GR5} (p. 384 No. 3.471.9)

\begin{eqnarray}
S_{1}^{\eta_{1}0\eta_{2}0}\left(\mathbf{k};0,\mathbf{x}_{2}\right)\  & = & \pi^{1/2}\int_{0}^{1}\,d\tau{\displaystyle \frac{1}{\tau{}^{1/2}}}e^{-i\mathbf{k}\cdot\mathbf{x}_{2}\tau}\int_{0}^{\infty}\,d\rho_{1}{\displaystyle \frac{1}{\rho_{1}^{\;1/2+1}}}\nonumber \\
 & \times & \exp\left(-\tau x_{2}^{2}\rho_{1}-\left(\frac{\tau\eta_{1}^{2}}{(4\tau)\rho_{1}}+\frac{(1-\tau)\left(k^{2}\tau+\eta_{2}^{2}\right)}{(4\tau)\rho_{1}}\right)\right)\nonumber \\
 & = & \pi^{1/2}\int_{0}^{1}\,d\tau e^{-i\mathbf{k}\cdot\mathbf{x}_{2}\tau}\frac{2\sqrt{\pi}\exp\left(-x_{2}\sqrt{(1-\tau)\left(k^{2}\tau+\eta_{2}^{2}\right)+\eta_{1}^{2}\tau}\right)}{\sqrt{(1-\tau)\left(k^{2}\tau+\eta_{2}^{2}\right)+\eta_{1}^{2}\tau}}\nonumber \\
 & = & 2\pi\int_{0}^{1}\,d\tau e^{-i\mathbf{k}\cdot\mathbf{x}_{2}\tau}\frac{\exp\left(-x_{2}L\right)}{L}\;,\label{eq:tau_int_L}
\end{eqnarray}
 where

\begin{equation}
L=\sqrt{(1-\tau)\left(k^{2}\tau+\eta_{2}^{2}\right)+\eta_{1}^{2}\tau}\:.\label{eq:L-1}
\end{equation}
 Interestingly enough, we see that when some more-complicated transition
amplitudes have been reduced to a single remaining integral, the integrand
can sometimes again take on a form akin to Slater orbitals, but with
no angular term in its argument and a complicated dependence on the
integration variable. We will refer to such as Yukawa-form functions.

Cheshire \cite{Cheshire} eq. (19) has setup, but not fully solved,
the related integral 
\begin{quotation}
\begin{eqnarray}
\hspace{-0.6cm}\left.I_{1}=\frac{\eta_{1}^{3/2}}{\sqrt{\pi}}S_{1s\,1}^{\eta_{1}0\eta_{2}0}\left(\frac{1}{2}\mathbf{k}_{f};0,\mathbf{x}_{2}\right)\right|_{\eta_{1}=\eta_{2}=1}\hspace{-0.3cm} & = & \left.\int d^{3}x_{1}\frac{\eta_{2}^{3/2}}{\sqrt{\pi}}e^{-\eta_{2}x_{12}}\frac{\eta_{1}^{3/2}}{\sqrt{\pi}}\frac{e^{-\eta_{1}x_{1}}}{x_{1}}e^{-i\frac{1}{2}\mathbf{k_{f}}\cdot\mathbf{x}_{1}}\right|_{\eta_{1}=\eta_{2}=1}\nonumber \\
 & = & \left.\frac{\eta_{1}^{3/2}}{\sqrt{\pi}}\frac{\eta_{2}^{3/2}}{\sqrt{\pi}}\left(-\frac{\partial}{\partial\eta_{2}}\right)\int d^{3}x_{1}\frac{e^{-\eta_{2}x_{12}}}{x_{12}}\frac{e^{-\eta_{1}x_{1}}}{x_{1}}e^{-i\frac{1}{2}\mathbf{k_{f}}\cdot\mathbf{x}_{1}}\right|_{\eta_{1}=\eta_{2}=1}\nonumber \\
 & = & \left.\frac{\eta_{1}^{3/2}}{\sqrt{\pi}}\frac{\eta_{2}^{3/2}}{\sqrt{\pi}}\left(-\frac{\partial}{\partial\eta_{2}}\right)S_{1}^{\eta_{1}0\eta_{2}0}\left(\frac{1}{2}\mathbf{k_{f}},;0,x_{2}\right)\right|_{\eta_{1}=\eta_{2}=1}\hspace{0.3cm}\,.\label{eq:cheshire}
\end{eqnarray}
\end{quotation}
Cheshire's 1964 \cite{Cheshire} Yukawa-form final integral has defied
reduction to analytic form for more than sixty years. The motivation
for the present paper was to show how to create an addition theorem
for this class of problems that allows the final integral to be done
term by term for $k\leq1$.

\section{A Different Yukawa-form Addition Theorem as a Generator of Macdonald
functions}

We wish to prove the following:

\textbf{Theorem 1.} For $k\leq1$,

\begin{equation}
\begin{array}{ccc}
{\displaystyle \frac{e^{-x_{2}\sqrt{Bk^{2}+C}}}{\sqrt{Bk^{2}+C}}} & = & {\displaystyle \sum_{n=0}^{\infty}}\sqrt{{\displaystyle \frac{2}{\pi}}}{\displaystyle \frac{\left(-1\right)^{n}B^{n}k^{2n}}{n!}}2^{-n}x_{2}^{n+1/2}C^{-n/2-1/4}K_{n+1/2}\left(x_{2}\sqrt{C}\right)\end{array}\label{eq:theorem1}
\end{equation}
 \textbf{Proof of Theorem 1.}

\begin{comment}
eq : yukBessekKseries is an inch down in the file ``new addition
theorem expand in C redo m7.nb''
\end{comment}

Consider the general Yukawa-form function,

\begin{eqnarray}
\frac{\exp\left(-x_{2}L\right)}{L}\:,\label{eq:cheshireKernel}
\end{eqnarray}
 where 
\begin{equation}
L\equiv\sqrt{k^{2}B+C}\quad\label{eq:L as B and C}
\end{equation}
 could, but need not, have the complicated dependence on the integration
variable $\tau$ in the argument seen in

\begin{equation}
L=\sqrt{k^{2}(1-\tau)\tau+\left(\eta_{1}^{2}-\eta_{2}^{2}\right)\tau+\eta_{2}^{2}}\quad.\label{eq:L as tau}
\end{equation}
 The argument in this example shares a problematic nature with (\ref{eq:yukseriesx1x2><})
in that it is non-integrable precisely because it is a quadratic form
confined within a square root.

An initial solution attempt using a Taylor series expansion in \emph{k}
gave terms that had echoes of Macdonald functions, but the coefficients
of the general term were difficult to discern. (This was one of those
interesting times when one's training in graphic design can lend service
to abstract mathematical visualization.) But given that any Yukawa-form
function is proportional to a Macdonald function with index 1/2, it
occurred to me that performing another Gaussian transform on it, \cite{GR5}
(p. 384 No. 3.471.9) splitting off and expanding the k-dependent term
in a series, \cite{GR5} (p. 27 No. 1.211.1) and performing the inverse
Gaussian transform on the remainder might lead to the correct series
of Macdonald functions with increasing indices.

Thus,

\begin{equation}
\begin{array}{ccc}
{\displaystyle \frac{e^{-x_{2}\sqrt{Bk^{2}+C}}}{\sqrt{Bk^{2}+C}}} & = & {\displaystyle \frac{1}{\sqrt{\pi}}}\int_{0}^{\infty}\,d\rho{\displaystyle \frac{e^{-\left(Bk^{2}+C\right)\rho}e^{-x_{2}^{2}/4/\rho}}{\rho^{\;1/2}}}\\
 & = & {\displaystyle \frac{1}{\sqrt{\pi}}}\int_{0}^{\infty}\,d\rho{\displaystyle \sum_{n=0}^{\infty}}{\displaystyle \frac{\left(-1\right)^{n}B^{n}k^{2n}\rho^{n+1/2-1}}{n!}}e^{-C\rho-x_{2}^{2}/4/\rho}\:,\\
 & = & {\displaystyle \sum_{n=0}^{\infty}}\sqrt{{\displaystyle \frac{2}{\pi}}}{\displaystyle \frac{\left(-1\right)^{n}B^{n}k^{2n}}{n!}}2^{-n}x_{2}^{n+1/2}C^{-n/2-1/4}K_{n+1/2}\left(x_{2}\sqrt{C}\right)
\end{array}\label{eq:yukBessekKseries}
\end{equation}
 which completes the proof.$\square$ For the arbitrary values $\left\{ C\to0.11,B\to0.13,x_{2}\to0.17,k\to0.23\right\} $
the first four terms in the series ${2.79367-0.051348+0.001318-0.000038}$
converge rapidly to the value on the left-hand side, $2.7436$.

\section{A Practical Application}

Let us now apply (\ref{eq:theorem1}) to the seemingly insoluble integral
(\ref{eq:tau_int_L}),

\begin{eqnarray}
S_{1}^{\eta_{1}0\eta_{2}0}\left(\mathbf{k};0,\mathbf{x}_{2}\right) & = & 2\pi\int_{0}^{1}\,d\tau e^{-i\mathbf{k}\cdot\mathbf{x}_{2}\tau}\frac{\exp\left(-x_{2}L\right)}{L}\label{eq:cheshireK}\\
 & = & 2\pi{\displaystyle \sum_{n=0}^{\infty}}{\displaystyle \frac{1}{\sqrt{\pi}}}{\displaystyle \frac{\left(-1\right)^{n}k^{2n}}{n!}}2^{-n+1/2}x_{2}^{n+1/2}\int_{0}^{1}\,d\tau e^{-i\mathbf{k}\cdot\mathbf{x}_{2}\tau}(1-\tau)^{n}\tau^{n}\nonumber \\
 & \times & \left(\left(\eta_{1}^{2}-\eta_{2}^{2}\right)\tau+\eta_{2}^{2}\right)^{-n/2-1/4}K_{n+1/2}\left(x_{2}\sqrt{\left(\eta_{1}^{2}-\eta_{2}^{2}\right)\tau+\eta_{2}^{2}}\right)\quad,\nonumber 
\end{eqnarray}
 For the first two terms, this can be directly integrated using the computer-calculus program \emph{Mathematica 7}, but for the remainder we change variables to $s=\sqrt{\left(\eta_{1}^{2}-\eta_{2}^{2}\right)\tau+\eta_{2}^{2}}$
to get

\begin{eqnarray}
S_{1}^{\eta_{1}0\eta_{2}0}\left(\mathbf{k};0,\mathbf{x}_{2}\right) & = & 2\pi{\displaystyle \sum_{n=0}^{\infty}}{\displaystyle \frac{1}{\sqrt{\pi}}}{\displaystyle \frac{\left(-1\right)^{n}k^{2n}}{n!}}2^{-n+1/2}x_{2}^{n+1/2}\label{eq:cheshireKs}\\
 & \times & \frac{2^{\frac{3}{2}-n}k^{2n}\left(\eta_{1}-\eta_{2}\right){}^{-2n}\left(\eta_{1}+\eta_{2}\right){}^{-2n}}{\sqrt{\pi}\left(\eta_{2}^{2}-\eta_{1}^{2}\right)n!}exp\left(-\frac{i\eta_{2}^{2}k\cdot x_{2}}{\eta_{2}^{2}-\eta_{1}^{2}}\right)\nonumber \\
 & \times & \int_{\eta_{2}}^{\eta_{1}}s^{\frac{1}{2}-n}x_{2}^{n+\frac{1}{2}}\left(s^{2}-\eta_{1}^{2}\right){}^{n}\left(s^{2}-\eta_{2}^{2}\right){}^{n}K_{n+\frac{1}{2}}\left(sx_{2}\right)exp\left(\frac{is^{2}k\cdot x_{2}}{\eta_{2}^{2}-\eta_{1}^{2}}\right)\,ds
\end{eqnarray}
 For $n=0$ we have

\begin{eqnarray}
S_{1}^{\eta_{1}0\eta_{2}0}\left(\mathbf{k};0,\mathbf{x}_{2}\right) & = & \frac{(1-i)\sqrt{2}\pi^{3/2}\exp\left(-\frac{ikx_{2}\eta_{2}^{2}}{\eta_{2}^{2}-\eta_{1}^{2}}-\frac{ix_{2}\left(\eta_{1}^{2}-\eta_{2}^{2}\right)}{4k}\right)}{\sqrt{\eta_{1}^{2}-\eta_{2}^{2}}\sqrt{k\cdot x_{2}}}\label{eq:cheshireKs0}\\
 & \times & \left(\text{erf}\left(\frac{(-1)^{3/4}\sqrt{x_{2}}\left(\eta_{1}^{2}+\eta_{2}\left(-\eta_{2}+2ik\right)\right)}{2\sqrt{k}\sqrt{\eta_{1}^{2}-\eta_{2}^{2}}}\right)\right.\nonumber \\
 & - & \left.\text{erf}\left(\frac{(-1)^{3/4}\sqrt{x_{2}}\left(2ik\eta_{1}+\eta_{1}^{2}-\eta_{2}^{2}\right)}{2\sqrt{k}\sqrt{\eta_{1}^{2}-\eta_{2}^{2}}}\right)\right)+\cdots.
\end{eqnarray}

\begin{comment}
This is the form given at \%114 in ``yukawa potential to a bessel
function 3 has transformation to sp m7.nb'' 
\end{comment}
When all of the parameters are of order one or less, this first term
accounts for 99\% of the value of the integral. Numerical integration
of the first line of (\ref{eq:cheshireK}) with $\left\{ \eta_{1}\to0.82,\eta_{2}\to0.66,x_{2}\to0.036,k\to0.019\right\} $,
for instance gives $6.4564-0.210837i$ and the $n=0$ term gives an
over-estimate of $6.50124-0.212271i$. Adding in the $n=1$ contribution
( $-0.0452324+0.00144522i$) gives an additional two decimal places
of accuracy: $6.45601-0.210825i$. %
\begin{comment}
Numeircal values from ``sp numerical check m7.nb''
\end{comment}
{} The $n=2$ and $3$ terms contribute $5\times10^{-4}$ and $-5\times10^{-6}$,
respectively.

More generally, we can write the half-integer Macdonald function as
a finite series

\begin{equation}
K_{n+\frac{1}{2}}\left(sx_{2}\right)\equiv\sqrt{\frac{\pi}{2}}\frac{e^{-sx_{2}}}{\sqrt{s}\sqrt{x_{2}}}\sum_{J=0}^{n}\frac{(J+n)!}{J!(n-J)!}\frac{1}{2^{J}s^{J}x_{2}^{J}}\:, \label{eq:half-integer Macdonald function}
\end{equation}
complete the square in the exponential $exp\left(-sx_{2}+\frac{is^{2}k\cdot x_{2}}{\eta_{2}^{2}-\eta_{1}^{2}}\right)$,
change variables to $s'=s+\frac{i\left(\eta_{2}^{2}-\eta_{1}^{2}\right)}{2k}$
with unit Jacobian, and expand the binomials to give an integrable
function, whose value is

\begin{eqnarray}
S_{1}^{\eta_{1}0\eta_{2}0}\left(\mathbf{k};0,\mathbf{x}_{2}\right) & = & \frac{2^{1-n}k^{2n}x_{2}^{n}\left(\eta_{1}-\eta_{2}\right){}^{-2n}\left(\eta_{1}+\eta_{2}\right){}^{-2n}e^{-\frac{i\eta_{2}^{2}k\cdot x_{2}}{\eta_{2}^{2}-\eta_{1}^{2}}}}{\left(\eta_{2}^{2}-\eta_{1}^{2}\right)n!}\sum_{n=0}^{\infty}\nonumber \\
 & \times & \sum_{m=0}^{n}(-1)^{m}\eta_{1}^{2m}\binom{n}{m}\sum_{j=0}^{n}(-1)^{j}\eta_{2}^{2j}\binom{n}{j}\sum_{J=0}^{n}\frac{2^{-J}x_{2}^{-J}(J+n)!e^{-\frac{ix_{2}\left(\eta_{1}^{2}-\eta_{2}^{2}\right)}{4k}}}{J!(n-J)!}\nonumber \\
 & \times & \sum_{K=0}^{-2j-J-2m+3n}\left(-\frac{1}{2}\right)^{K}\left(\frac{i\left(\eta_{2}^{2}-\eta_{1}^{2}\right)}{k}\right){}^{K}\binom{-2j-J-2m+3n}{K}\nonumber \\
 & \times & \left(\frac{\left(\frac{1}{\eta_{1}^{2}-\eta_{2}^{2}}\right){}^{\frac{1}{2}(2j+J+K+2m-3n-1)}\left(\eta_{1}^{2}+\eta_{2}\left(-\eta_{2}+2ik\right)\right){}^{2j+J+K+2m-3n}}{2kx_{2}\left(\eta_{2}\left(2k+i\eta_{2}\right)-i\eta_{1}^{2}\right){}^{2j+J+K+2m-3n}}\right.\nonumber \\
 & \times & \left(-ix_{2}\right){}^{\frac{1}{2}(2j+J+K+2m-3n+1)}\nonumber \\
 & \times & \Gamma\left(\frac{1}{2}(-2j-J-K-2m+3n+1),-\frac{ix_{2}\left(\eta_{1}^{2}+\left(2ik-\eta_{2}\right)\eta_{2}\right){}^{2}}{4k\left(\eta_{1}^{2}-\eta_{2}^{2}\right)}\right)\nonumber \\
 & - & \frac{1}{2}\left(\frac{2k\eta_{1}+i\left(\eta_{2}^{2}-\eta_{1}^{2}\right)}{k}\right)^{-2j-J-K-2m+3n+1}\nonumber \\
 & \times & \left.\Gamma\left(\frac{1}{2}(-2j-J-K-2m+3n-1)+1,\frac{ix_{2}\left(2k\eta_{1}+i\left(\eta_{2}^{2}-\eta_{1}^{2}\right)\right){}^{2}}{4k\left(\eta_{1}^{2}-\eta_{2}^{2}\right)}\right)\right)\:,
\end{eqnarray}

\noindent where the final term results from converting an Exponential
Integral function $E_{n}(z)$ %from the file sp numerical check reredo end m7.nb
to a confluent hypergeometric $\,_{1}F_{1}$ function \cite{functions.wolfram.com/06.34.26.0002.01}
and thence to an incomplete Gamma function, \cite{functions.wolfram.com/06.06.26.0002.01}
an alternative to the Error function \cite{functions.wolfram.com/06.25.26.0001.01}
of (\ref{eq:cheshireKs0}), above.

Returning to the Cheshire integral, which has $\eta_{2}\to\eta_{1}$,
the complexity of the above reduces considerably to

\begin{eqnarray}
S_{1}^{\eta_{1}0\eta_{1}0}\left(\mathbf{k};0,\mathbf{x}_{2}\right) & = & 2\pi{\displaystyle \sum_{n=0}^{\infty}}{\displaystyle \frac{1}{\sqrt{\pi}}}{\displaystyle \frac{\left(-1\right)^{n}k^{2n}}{n!}}2^{-n+1/2}x_{2}^{n+1/2}\eta_{1}^{-n-\frac{1}{2}}K_{n+\frac{1}{2}}\left(x_{2}\eta_{1}\right)\quad,\label{eq:cheshireK-1}\\
 & \times & \int_{0}^{1}\,d\tau e^{-i\mathbf{k}\cdot\mathbf{x}_{2}\tau}(1-\tau)^{n}\tau^{n}\nonumber \\
 & = & 2\pi\sum_{n=0}^{\infty}\frac{(-1)^{n}2^{-3n-\frac{1}{2}}k^{2n}x_{2}^{n+\frac{1}{2}}\eta_{1}^{-n-\frac{1}{2}}}{\Gamma\left(n+\frac{3}{2}\right)}\nonumber \\
 & \times & K_{n+\frac{1}{2}}\left(x_{2}\eta_{1}\right)\,_{1}F_{1}\left(n+1;2n+2;-ik\cdot x_{2}\right)\nonumber 
\end{eqnarray}
 We, thus, have developed a series solution to the Cheshire integral,
and its generalization, for $k\leq1$.

\section{A one-range Addition Theorem for Slater Orbitals}

The utility of this one-range addition theorem goes far beyond the
application, above, that motivated it. One can consider an application
of this new addition theorem, with k=1, for a conventional Slater orbital with a shifted-position
vector magnitude in spherical coordinates (wherein one has the freedome to choose the direction of the rotational axis such that  $\mathbf{x}_{1}\cdot\mathbf{x}_{2}=x_{1}x_{2}\cos\left(\theta\right)$),
grouped as

\textbf{Corollary 1.}

\begin{equation}
\begin{array}{ccc}
\left.{\displaystyle \frac{e^{-\eta\sqrt{\left(x_{1}^{2}-2\mathbf{x}_{1}\cdot\mathbf{x}_{2}\right)k^{2}+x_{2}^{2}}}}{\sqrt{\left(x_{1}^{2}+2\mathbf{x}_{1}\cdot\mathbf{x}_{2}\right)k^{2}+x_{2}^{2}}}}\right|_{k=1} & = & {\displaystyle \sum_{n=0}^{\infty}}{\displaystyle \frac{1}{\sqrt{\pi}}}{\displaystyle \frac{\left(-1\right)^{n}\left(x_{1}^{2}-2\mathbf{x}_{1}\cdot\mathbf{x}_{2}\right)^{n}}{n!}}2^{-n+1/2}\eta^{n+1/2}x_{2}^{-n-1/2}K_{n+1/2}\left(\eta x_{2}\right)\\
 & = & {\displaystyle \sum_{n=0}^{\infty}}{\displaystyle \frac{1}{\sqrt{\pi}}}{\displaystyle \frac{\left(-1\right)^{n}}{n!}}2^{-n+1/2}\eta^{n+1/2}x_{2}^{-n-1/2}K_{n+1/2}\left(\eta x_{2}\right)\\
 & \times & {\displaystyle \sum_{j=0}^{n}}\left(-1\right)^{j}2^{j}x_{2}^{j}\binom{n}{j}x_{1}^{2n-j}{\displaystyle \sum_{m=j,j-2,\ldots}\frac{(2m+1)j!2^{\frac{m-j}{2}}}{\frac{j-m}{2}!(j+m+1)\text{!!}}}P_{m}(\cos(\theta))
\end{array}\;,\label{eq:yukseries x1^2-2ux1x2}
\end{equation}
by setting$\left\{ x_{2}\to\eta,C\to x_{2}^{2},B\to x_{1}^{2}-2\mathbf{x}_{1}\cdot\mathbf{x}_{2},\,k\to1\right\} $
in Theorem 1 (\ref{eq:theorem1}). Note that this one-range addition
theorem has a finite second series rather than an infinite one such as that 
in (\ref{eq:one-range}) (that has additional finite series, too).
One sees that the expansion of powers of $\cos\left(\theta\right)$
in Legendre Polynomials \textendash{} given in Weisstein \cite{mathworld.wolfram.com/LegendrePolynomial.html}
(where the smallest value of \emph{m} is zero or one depending on
whether \emph{j} is even or odd, resp.) \textendash{} makes this version
seem less appealing than the one found in Magnus, Oberhettinger, and
Soni \cite{Magnus Oberhettinger and Soni}, (\ref{eq:yukseriesx1x2><}),
above, since applying orthogonality relations will truncate the finite
series in \emph{m} rather than the infinite series in \emph{n}. But
in many cases, it will avoid the significant negative numerical consequences
of round-off errors in cancelling terms of large magnitude and opposite
signs that sometimes appear in applications of the latter two-range
addition theorem.

It is interesting that (\ref{eq:yukseries x1^2-2ux1x2}) is not the
only grouping possible, so that we essential have four one-range addition
theorems in one for Slater orbitals in spherical coordinates. One
may set $\left\{ x_{2}\to\eta,\,C\to x_{1}^{2}\right.,$ $\left.B\to x_{2}^{2}-2\mathbf{x}_{1}\cdot\mathbf{x}_{2},k\to1\right\} $
in Theorem 1 (\ref{eq:theorem1}) to obtain

\textbf{Corollary 2.}

\begin{equation}
\begin{array}{ccc}
\left.{\displaystyle \frac{e^{-\eta\sqrt{x_{1}^{2}+\left(-2\mathbf{x}_{1}\cdot\mathbf{x}_{2}+x_{2}^{2}\right)k^{2}}}}{\sqrt{x_{1}^{2}+\left(-2\mathbf{x}_{1}\cdot\mathbf{x}_{2}+x_{2}^{2}\right)k^{2}}}}\right|_{k=1} & = & {\displaystyle \sum_{n=0}^{\infty}{\displaystyle \frac{1}{\sqrt{\pi}}}\frac{(-1)^{n}\left(-2\mathbf{x}_{1}\cdot\mathbf{x}_{2}+x_{2}^{2}\right)^{n}}{n!}2^{\frac{1}{2}-n}x_{1}^{-n-\frac{1}{2}}\eta^{n+\frac{1}{2}}K_{n+\frac{1}{2}}\left(x_{1}\eta\right)}\end{array}\label{eq:yukseries x2^2-2ux1x2}
\end{equation}
Since the Macdonald function on the right-hand side has half-integer
indices, it is an exponential multiplying inverse powers of its argument.
So this addition theorem essentially allows one to pluck out problematic
pieces of the quadratic from in the square root of the exponential
leaving an integrable exponential (if we are integrating over $x_{1}$).
The cost of dealing with an infinite series may be worth it if $x_{2}$
is small enough that the series converges reasonably fast.

Or one might even decide to pluck out everything \emph{except} the
problematic angular portion of the quadratic from in the square root
of the exponential by setting $\left\{ x_{2}\to\eta,\,C\to-2\mathbf{x}_{1}\cdot\mathbf{x}_{2},B\to x_{1}^{2}+x_{2}^{2},\right.$
$\left.k\to1\right\} $ to obtain

\textbf{Corollary 3.}

\begin{equation}
\left.{\displaystyle \frac{e^{-\eta\sqrt{\left(x_{1}^{2}+x_{2}^{2}\right)k^{2}-2\mathbf{x}_{1}\cdot\mathbf{x}_{2}}}}{\sqrt{\left(x_{1}^{2}++x_{2}^{2}\right)k^{2}-2\mathbf{x}_{1}\cdot\mathbf{x}_{2}}}}\right|_{k=1}\hspace{-0.3cm}=\hspace{-0cm}{\displaystyle \sum_{n=0}^{\infty}{\displaystyle \frac{1}{\sqrt{\pi}}}\frac{(-1)^{n}\left(x_{1}^{2}+x_{2}^{2}\right)^{n}}{n!}}2^{\frac{1}{4}-\frac{3n}{2}}\eta_{2}^{n+\frac{1}{2}}\left(-\mathbf{x}_{1}\cdot\mathbf{x}_{2}\right){}^{-\frac{n}{2}-\frac{1}{4}}K_{n+\frac{1}{2}}\left(\sqrt{-2\mathbf{x}_{1}\cdot\mathbf{x}_{2}}\eta_{2}\right)\label{eq:yukseries x2^2-2ux1x2-1}
\end{equation}
While this choice might seem too strange \textendash{} with its imaginary
argument in the Macdonald function \textendash{} to pursue further,
it has secondary value in that one can actually integrate the $n=0$
term over the angular variable $u=\cos\left(\theta\right),$

\textbf{Theorem 2.}

\begin{equation}
\int_{-1}^{1}\frac{e^{-\sqrt{2}\eta_{2}\sqrt{-ux_{1}x_{2}}}}{\sqrt{-ux_{1}x_{2}}}\,du=\frac{\sqrt{2}\left(-e^{-\sqrt{2}\sqrt{x_{1}x_{2}}\eta_{2}}+e^{-i\sqrt{2}\sqrt{x_{1}x_{2}}\eta_{2}}\right)}{x_{1}x_{2}\eta_{2}}\:,\label{eq:n=00003D00003D0 angular int}
\end{equation}
using  \emph{Mathematica 7}. This adds
to our set of angular integrals that we have not found in the literature
\textendash{} appearing as equation (24) in a previous paper.\cite{Casey}

Finally, one can take the complement of Corollary 3 by setting $\left\{ x_{2}\to\eta,\,C\to x_{1}^{2}+x_{2}^{2},B\to-2\mathbf{x}_{1}\cdot\mathbf{x}_{2},\right.$
$\left.k\to1\right\} $ to obtain

\textbf{Corollary 4.}

\begin{equation}
\left.{\displaystyle \frac{e^{-\eta\sqrt{x_{1}^{2}+x_{2}^{2}+\left(-2\mathbf{x}_{1}\cdot\mathbf{x}_{2}\right)k^{2}}}}{\sqrt{x_{1}^{2}+x_{2}^{2}+\left(-2\mathbf{x}_{1}\cdot\mathbf{x}_{2}\right)k^{2}}}}\right|_{k=1}\hspace{-0.2cm}=\hspace{-0.1cm}{\displaystyle \sum_{n=0}^{\infty}{\displaystyle \frac{1}{\sqrt{\pi}}}\frac{(-1)^{n}\left(-\mathbf{x}_{1}\cdot\mathbf{x}_{2}\right){}^{n}}{n!}}\sqrt{2}\left(x_{1}^{2}+x_{2}^{2}\right){}^{\frac{1}{2}\left(-n-\frac{1}{2}\right)}\eta_{2}^{n+\frac{1}{2}}K_{n+\frac{1}{2}}\left(\sqrt{x_{1}^{2}+x_{2}^{2}}\eta\right)\label{eq:corollary4}
\end{equation}
Because this series involving the half-integer the MacDonald function
has terms of order $x_{j}^{-n-1}$ to $x_{j}^{-1}$ (multiplying a
negative exponential of both variables), it should be the most reliably
convergent of these four one-range addition theorems. For this reason
we will use this one-range addition theorem as a test-case to solve
the integral over two Slater orbitals (\ref{eq:theorem3}) in the
next section.

\section{A Test of a this One-range Addition Theorem for Slater Orbitals}

Let us use this one-range addition theorem (\ref{eq:corollary4})
in the second product within the integral over two Slater orbitals
(\ref{eq:SVxVxyG}). The result is a doubly infinite series whose
sum is the known analytic function,

\textbf{Theorem 3.}

\begin{eqnarray}
S_{1}^{\eta_{1}0\eta_{2}0}\left(0;0,\mathbf{x}_{2}\right) & = & \frac{4\pi\left(e^{-x_{2}\eta_{2}}-e^{-x_{2}\eta_{1}}\right)}{x_{2}\left(\eta_{1}^{2}-\eta_{2}^{2}\right)}\nonumber \\
 & = & \sum_{n=0}^{\infty}\sum_{k=0}^{\infty}\left(\sum_{j=0}^{\left\lfloor \frac{|n-1|}{2}-\frac{1}{2}\right\rfloor }\sum_{i=0}^{\frac{n}{2}}\frac{\sqrt{\pi}i^{3n}\eta_{1}(-1)^{k-i}2^{-j+\frac{n}{2}+3}\Gamma\left(\frac{n+3}{2}\right)\binom{\frac{n}{2}}{i}\eta_{2}^{n-2i}\left(\frac{n+3}{2}\right)_{k}}{j!k!\Gamma(n+2)\left(\frac{1}{2}(|n-1|-2j-1)\right)!}\right.\nonumber \\
 & \times & \left(\frac{1}{2}(|n-1|+2j-1)\right)!\nonumber \\
 & \times & \left.x_{2}^{-2i+\frac{1}{2}(n-2j)+j+2k+\frac{n}{2}+2}\left(\eta_{1}^{2}-\eta_{2}^{2}\right){}^{k}\Gamma\left(2i-j-2k-\frac{n}{2}-2,x_{2}\eta_{2}\right)\right)\;,\label{eq:theorem3}
\end{eqnarray}

\noindent
where the \emph{n}-sum is over even values only so that the  \emph{floor} function that is the upper limit of the \emph{j}-sum, $\left\lfloor \frac{| n-1| }{2}-\frac{1}{2}\right\rfloor$, gives $0$ for $n=0$ and $n/2-1$ for $n>=2$.

\noindent
 \textbf{Proof of Theorem 3.}

Since there is no angular dependence in the first Slater orbital,

\begin{eqnarray}
S_{1}^{\eta_{1}0\eta_{2}0}\left(0;0,\mathbf{x}_{2}\right) & = & \int d^{3}x_{1}\frac{e^{-\eta_{1}x_{1}}}{x_{1}}\frac{e^{-\eta_{2}x_{12}}}{x_{12}}\nonumber \\
 & = & 2\pi\int_{0}^{\infty}dx_{1}x_{1}^{2}\frac{e^{-\eta_{1}x_{1}}}{x_{1}}\int_{-1}^{1}du{\displaystyle \sum_{n=0}^{\infty}}{{\displaystyle \frac{1}{\sqrt{\pi}}}\frac{(-1)^{n}\left(-ux_{1}x_{2}\right){}^{n}}{n!}}\nonumber \\
 & \times & \hspace{1.6cm}\sqrt{2}\left(x_{1}^{2}+x_{2}^{2}\right){}^{\frac{1}{2}\left(-n-\frac{1}{2}\right)}\eta_{2}^{n+\frac{1}{2}}K_{n+\frac{1}{2}}\left(\sqrt{x_{1}^{2}+x_{2}^{2}}\eta_{2}\right)\nonumber \\
 & = & 2\pi\int_{0}^{\infty}dx_{1}x_{1}^{2}\frac{e^{-\eta_{1}x_{1}}}{x_{1}}{\displaystyle \sum_{n=0}^{\infty}}{{\displaystyle \frac{1}{\sqrt{\pi}}}\frac{(-1)^{n}}{n!}\frac{\left((-1)^{n}+1\right)x_{1}^{n}x_{2}^{n}}{n+1}}\nonumber \\
 & \times & \hspace{1.6cm}\sqrt{2}\left(x_{1}^{2}+x_{2}^{2}\right){}^{\frac{1}{2}\left(-n-\frac{1}{2}\right)}\eta_{2}^{n+\frac{1}{2}}K_{n+\frac{1}{2}}\left(\sqrt{x_{1}^{2}+x_{2}^{2}}\eta_{2}\right)\:,\label{eq:Sone-range}
\end{eqnarray}

\noindent where the factor $\left((-1)^{n}+1\right)$ on the fourth
line gives a factor of two for even \emph{n} and zero for odd\emph{
n} so that the $(-1)^{n}$ that precedes it is universally \textendash{}
thus redundantly \textendash{} one. Since the Gaussian transform of
the Macdonald function is known \cite{PBM5}(p. 230 No. 3.16.1.13),

\begin{eqnarray}
p^{-\mu}K_{\nu}\left(a\sqrt{p}\right) & = & \int_{0}^{\infty}2^{-\nu-1}a^{\nu}e^{-\frac{a^{2}}{4\rho_{2}}-p\rho_{2}}\rho_{2}^{\mu-\frac{\nu}{2}-1}U\left(\mu+\frac{\nu}{2},\nu+1,\frac{a^{2}}{4\rho_{2}}\right)\,d\rho_{2}\:,\label{eq:Gaussian to Madonald}\\
 &  & \left\{ Re\:p>0,\left|arg\,a\right|<\pi/4\:\mathrm{or\:\left|arg\,a\right|=\pi/4,}\:Re\,\mu>-1/2\right\} \nonumber 
\end{eqnarray}

\noindent where U is the Tricomi confluent hypergeometric function,
one might as well use that approach in the $x_{1}$ integral. The
Gaussian transform of $\frac{e^{-\eta_{1}x_{1}}}{x_{1}}$ is given
in (\ref{eq:squarecompleted}).

Since there is no cross term and no momentum vector k in the present
integral, the equivalent of the exponential in the last line of (\ref{eq:squarecompleted})
is $exp\left(-\left(\rho_{1}+\rho_{2}\right)x_{1}^{2}-x_{2}^{2}\rho_{2}\right)$
and no change of variables is required to integrate over $x_{1}$:

\begin{eqnarray}
S_{1}^{\eta_{1}0\eta_{2}0}\left(0;0,\mathbf{x}_{2}\right) & \hspace{-0.2cm}= & \hspace{-0.2cm}\sum_{n=0}^{\infty}\int_{0}^{\infty}\int_{0}^{\infty}d\rho_{2}d\rho_{1}\frac{\left((-1)^{n}+1\right)x_{1}^{n}2^{n+1}x_{2}^{n}\rho_{2}^{n-\frac{1}{2}}\left(\rho_{1}+\rho_{2}\right){}^{-\frac{n}{2}-\frac{3}{2}}\Gamma\left(\frac{n+3}{2}\right)e^{-\rho_{2}x_{2}^{2}-\frac{\eta_{1}^{2}}{4\rho_{1}}-\frac{\eta_{2}^{2}}{4\rho_{2}}}}{\sqrt{\rho_{1}}\Gamma(n+2)}\nonumber \\
 & = & \sum_{n=0}^{\infty}\int_{0}^{1}d\tau_{1}\int_{0}^{\infty}d\rho_{2}\frac{\left((-1)^{n}+1\right)x_{1}^{n}2^{n+1}x_{2}^{n}\rho_{2}^{\frac{n}{2}-\frac{3}{2}}\left(1-\tau_{1}\right){}^{n/2}\Gamma\left(\frac{n+3}{2}\right)}{\sqrt{\tau_{1}}\Gamma(n+2)}\nonumber \\
 & \times & \hspace{4.5cm}\exp\left(-x_{2}^{2}\rho_{2}-\left(\frac{\eta_{1}^{2}\left(1-\tau_{1}\right)}{4\tau_{1}}+\frac{\eta_{2}^{2}}{4}\right)/\rho_{2}\right)\:,\label{eq:Sone-rangetau pho2}
\end{eqnarray}

\noindent where we have used the same change of variables to $\tau_{1}$
as in (\ref{eq:tau transform}). The $\rho_{2}$ integral is easily
done \cite{GR5}(p. 384 No. 3.471.9), giving

\begin{eqnarray}
S_{1}^{\eta_{1}0\eta_{2}0}0;0,\mathbf{x}_{2} & = & \sum_{n=0}^{\infty}\int_{0}^{1}d\tau_{1}\frac{\left((-1)^{n}+1\right)x_{1}^{n}2^{\frac{n}{2}+\frac{5}{2}}x_{2}^{\frac{n}{2}+\frac{1}{2}}\left(1-\tau_{1}\right){}^{n/2}\Gamma\left(\frac{n+3}{2}\right)}{\sqrt{\tau_{1}}\Gamma(n+2)}\left(\frac{\eta_{2}^{2}\tau_{1}-\eta_{1}^{2}\left(\tau_{1}-1\right)}{\tau_{1}}\right){}^{\frac{n-1}{4}}\nonumber \\
 & \times & \hspace{2.5cm}K_{\frac{1-n}{2}}\left(x_{2}\sqrt{\frac{\eta_{2}^{2}\tau_{1}-\eta_{1}^{2}\left(\tau_{1}-1\right)}{\tau_{1}}}\right)\nonumber \\
 & = & \sum_{n=0}^{\infty}\int_{\eta_{2}}^{\infty}\,ds\frac{\left((-1)^{n}+1\right)x_{1}^{n}2^{\frac{n}{2}+\frac{7}{2}}\eta_{1}s^{\frac{n}{2}+\frac{1}{2}}x_{2}^{\frac{n}{2}+\frac{1}{2}}\Gamma\left(\frac{n+3}{2}\right)\left(s^{2}-\eta_{2}^{2}\right)^{n/2}\left(s^{2}+\eta_{1}^{2}-\eta_{2}^{2}\right){}^{-\frac{n}{2}-\frac{3}{2}}}{\Gamma(n+2)}\nonumber \\
 & \times & \hspace{2.5cm}K_{\frac{1-n}{2}}\left(sx_{2}\right)\;,\label{eq:Stau}
\end{eqnarray}

Here we changed variables to a somewhat different $s=\frac{\sqrt{\eta_{2}^{2}\tau_{1}-\eta_{1}^{2}\left(\tau_{1}-1\right)}}{\sqrt{\tau_{1}}}$
for this new problem than in (\ref{eq:cheshireKs}), but to the same
purpose. This integral can be done if we expand $\left(s^{2}-\eta_{2}^{2}\right)^{n/2}$
into a finite sum of isolate powers of \emph{s}, since \emph{n} is
even, and expand the denominator $\left(s^{2}+\eta_{1}^{2}-\eta_{2}^{2}\right){}^{-\frac{n}{2}-\frac{3}{2}}$
into a series via \cite{PBM3} (p. 455 No. 7.3.1.27)

\begin{equation}
z^{a}=\,_{2}F_{1}(-a,b;b;1-z)\:.\label{eq:denom as 2F1}
\end{equation}
We would like, then, to explicate these powers of \emph{z} via \cite{PBM3}
(p. 430 No. 7.2.1.1)

\begin{eqnarray}
\,_{2}F_{1}(a,b;c;1-z) & = & \sum_{k=0}^{\infty}\frac{(1-z)^{k}(a)_{k}(b)_{k}}{k!(c)_{k}}\:,\label{eq:2F1 as pochh}\\
 &  & \text{\ensuremath{\left|1-z\right|}<1\ensuremath{\lor}(|1-z|=1\ensuremath{\land\Re}(-a-b+c)>0\ensuremath{\lor}(|1-z|=1\ensuremath{\land}1-z\ensuremath{\neq}1\ensuremath{\land}-1<\ensuremath{\Re}(-a-b+c)\ensuremath{\leq}0)|}
\end{eqnarray}
but because of the condition $|1-z|<1$ in this second step, we must
first rewrite

\begin{equation}
\left(s^{2}+\eta_{1}^{2}-\eta_{2}^{2}\right){}^{-\frac{n}{2}-\frac{3}{2}}=s^{-n-3}\left(1+\frac{\eta_{1}^{2}-\eta_{2}^{2}}{s^{2}}\right)^{-\frac{n}{2}-\frac{3}{2}}\:,\label{eq:get z<1}
\end{equation}

\noindent and transform the second factor. Since $\eta_{2}$ is typically
the nuclear charge, $s$ will be larger than one throughout the integral
and therefore \emph{$|1-z|=\left|-\frac{\eta_{1}^{2}-\eta_{2}^{2}}{s^{2}}\right|<1$}.
Finally, we expand $K_{\frac{1-n}{2}}\left(sx_{2}\right)$ into a
second finite series before integrating to obtain the desired result
(\ref{eq:theorem3}).$\square$

The analytic function in the first line of Theorem 3 (\ref{eq:theorem3})
has value $51.3025821$ for parameters arbitrarily chosen to be $\left\{ \eta_{1}\to0.11,\eta_{2}\to0.13,x_{2}\to0.17\right\} $.
For the $n=0$ term in the outermost series, we find that the first
three terms in the $k$ series contribute strongly: $\{39.1836,8.45916,2.008,0.499656,0.127812,0.033291,$
$0.008782,0.002339,0.000628\}$, and sum to $50.3232$. For $n=2$,
the terms in the $k$ series fall off more rapidly, $\{0.527506,0.082122,0.017650,0.004191,0.001043,0.000267,0.00007,0.000018,5\times10^{-6}\}$,
and sum to $0.632872$. For $n=4$ and higher, each term in the $k$
series again shows convergence at roughly a digit of accuracy every
other term, $\{0.115654,0.017087,0.003642,0.000862,0.000214,0.00005,0.000014,3\times10^{-6},1\times10^{-6}\}$,
in this case summing to $0.137535$. The next two terms in the $n$
series $0.0593952$ and $0.033087$ indicate fairly slow convergence
(to $51.1861$ with these five terms), but there are, indeed, no alternating
signs that might cause roundoff errors such as I found when using
the two-range addition theorem in Magnus, Oberhettinger, and Soni
\cite{Magnus Oberhettinger and Soni}, (\ref{eq:yukseriesx1x2><}),
applied to five Slater orbitals.%
\begin{comment}
Numbers from the file ``no j sum via gaussian general n c m7.nb''
\end{comment}

%numbers from the file "no j sum via gaussian general n c m7.nb"

This sequence of steps to obtain a doubly-infinite series approximation
to an analytic function obtainable exactly by other means such as
(\ref{eq:SVxVxyG}) may sound like an irrational method of working.
But those other means likely have reached their limit of applicability
with Fromm and Hill's \cite{Fromm and Hill} tour de force Fourier-transform-integration
over the angular and radial variables for a product of six Slater
orbitals in three 3D integration variables, three orbitals of which
had shifted coordinates. So what are we to do if we have seven Slater
orbitals in four 3D integration variables and/or four orbitals having
shifted coordinates? It could be that this one-range addition theorem
might allow us to trudge a viable path.

\subsection{The case of $\eta_{1}=\eta_{2}$}

There will be many applications for which the atomic charge (or screened
charges or nuclear decay lengths in the various applications noted
in the introduction) are the same for both Slater orbitals. One need
not entirely rederive (\ref{eq:SVxVxyG}) but may simply take the
limit of the final result, the first line of the following:

\textbf{Theorem 4.}

\begin{eqnarray}
S_{1}^{\eta_{2}0\eta_{2}0}\left(0;0,\mathbf{x}_{2}\right) & = & \frac{2\pi e^{-x_{2}\eta_{2}}}{\eta_{2}}\nonumber \\
 & =\sum_{n=0}^{\infty} & \sum_{j=0}^{\left\lfloor \frac{|n-1|}{2}-\frac{1}{2}\right\rfloor }\sum_{i=0}^{\frac{n}{2}}\left(\frac{\sqrt{\pi}(-1)^{-i}i^{3n}2^{-j+\frac{n}{2}+3}\Gamma\left(\frac{n+3}{2}\right)\binom{\frac{n}{2}}{i}\left(\frac{1}{2}(|n-1|+2j-1)\right)!}{j!\Gamma(n+2)\left(\frac{1}{2}(|n-1|-2j-1)\right)!}\right.\nonumber \\
 & \times & \left.x_{2}^{-2i+n+2}\eta_{2}^{-2i+n+1}\Gamma\left(2i-j-\frac{n}{2}-2,x_{2}\eta_{2}\right)\right)\;.\label{eq:etas equal}
\end{eqnarray}
 \textbf{Proof of Theorem 4.}

For the general case where $\eta_{1}$ is not necessarily equal to
$\eta_{2}$, to perform the final integral we had expand the denominator
$\left(1+\frac{\eta_{1}^{2}-\eta_{2}^{2}}{s^{2}}\right)^{-\frac{n}{2}-\frac{3}{2}}$
into a an infinite hypergeometric series (\ref{eq:denom as 2F1}).
But this step becomes unnecessary if $\eta_{1}=\eta_{2}$ since this
expression then becomes one. As before, we expand $K_{\frac{1-n}{2}}\left(sx_{2}\right)$
into a finite series before integrating to obtain the desired result
(\ref{eq:etas equal}).$\square$ The analytic function in the first
line of Theorem 4 (\ref{eq:etas equal}) has value $47.27577$ for
parameters chosen to be $\left\{ \eta_{1}\to0.13,\eta_{2}\to0.13,x_{2}\to0.17\right\} $,
whereas the first four terms in the series $\{46.3079,0.623416$,
$0.136682,0.0591038\}$ (summing to $47.1271$) show convergence at
roughly a digit every other term. %numbers from the file "no j sum via gaussian general n equal etas m7.nb"

The reader wishing to take such an approach with a different problem
is cautioned that if the expansions of $\left(s^{2}-\eta_{2}^{2}\right)^{n/2}$
and $K_{\frac{1-n}{2}}\left(sx_{2}\right)$ into finite sums are not
done before the integral is attempted, the result may very well include
almost-cancelling terms of infinite magnitude along with finite terms.

\section{This One-range Addition Theorem may also be used for Cartesian Coordinates}

There may be some problems for which Cartesian coordinates would be
more useful than spherical ones, in which case Corollary 3 becomes
(with the choice of coordinates aligning the z axis with the direction
of the shifted charge center), with $\left\{ x_{2}\to\eta,\,C\to\left(z_{1}-z_{2}\right){}^{2},\,B\to x_{1}^{2}+y_{1}^{2},\,k\to1\right\} $,

\textbf{Corollary 5.}

\begin{equation}
\left.\frac{e^{-\eta\sqrt{\left(x_{1}^{2}+y_{1}^{2}\right)k^{2}+\left(z_{1}-z_{2}\right){}^{2}}}}{\sqrt{\left(x_{1}^{2}+y_{1}^{2}\right)k^{2}+\left(z_{1}-z_{2}\right){}^{2}}}\right|_{k=1}\hspace{-0.3cm}=\hspace{-0cm}{\displaystyle \sum_{n=0}^{\infty}{\displaystyle \frac{1}{\sqrt{\pi}}}\frac{(-1)^{n}\left(x_{1}^{2}+y_{1}^{2}\right)^{n}}{n! 2^{n-\frac{1}{2}}}     K_{n+\frac{1}{2}}\left(\sqrt{\left(z_{1}-z_{2}\right){}^{2}}\eta\right)  \hspace{-0.1cm}  \left(\left(z_{1}-z_{2}\right){}^{2}\right){}^{\frac{1}{2}\left(-n-\frac{1}{2}\right)}\eta_{2}^{n+\frac{1}{2}},}\label{eq:corollary5}
\end{equation}
in which one must be careful to substitute $\sqrt{\left(z_{1}-z_{2}\right){}^{2}}=\left|z_{1}-z_{2}\right|$
since $z_{2}$ will be larger than $z_{1}$ at some points in an integral
over the latter. Unfortunately, this reverts to a two-range form.

The Cartesian-coordinate equivalent of Corollary 4, with $\left\{ x_{2}\to\eta,\,C\to x_{1}^{2}+y_{1}^{2},\,B\to\left(z_{1}-z_{2}\right){}^{2},\,k\to1\right\} $,
is

\textbf{Corollary 6.}

\begin{equation}
\left.\frac{e^{-\eta\sqrt{x_{1}^{2}+y_{1}^{2}+\left(z_{1}-z_{2}\right){}^{2}k^{2}}}}{\sqrt{x_{1}^{2}+y_{1}^{2}+\left(z_{1}-z_{2}\right){}^{2}k^{2}}}\right|_{k=1}\hspace{-0.2cm}=\hspace{-0.1cm}{\displaystyle \sum_{n=0}^{\infty}{\displaystyle \frac{1}{\sqrt{\pi}}}\frac{(-1)^{n}\left(z_{1}-z_{2}\right){}^{2n}}{n!}}2^{\frac{1}{2}-n}K_{n+\frac{1}{2}}\left(\sqrt{x_{1}^{2}+y_{1}^{2}}\eta\right)\left(x_{1}^{2}+y_{1}^{2}\right){}^{\frac{1}{2}\left(-n-\frac{1}{2}\right)}\eta_{2}^{n+\frac{1}{2}}.\label{eq:corollary6}
\end{equation}
The evaluation of this series follows that for its spherical-coordinate
equivalent, Corollary 4, through the substitution (\ref{eq:Gaussian to Madonald}).
But in this Cartesian-coordinate version, the $n=0$ term can be reduced
to analytic form. The $n=0$ term of (\ref{eq:corollary6}) becomes
\begin{equation}
\frac{\exp\left(\rho_{1}\left(-x_{1}^{2}-y_{1}^{2}-z_{1}^{2}\right)-\rho_{2}\left(x_{1}^{2}+y_{1}^{2}\right)-\frac{\eta_{1}^{2}}{4\rho_{1}}-\frac{\eta_{2}^{2}}{4\rho_{2}}\right)}{\pi\sqrt{\rho_{1}}\sqrt{\rho_{2}}}\:,\label{eq:n=00003D0cart}
\end{equation}
which can be simply integrated over $\left\{ z_{1},-\infty,\infty\right\} $,
$\left\{ y_{1},-\infty,\infty\right\} $, $\left\{ x_{1},-\infty,\infty\right\} $,
$\left\{ \rho_{2},0,\infty\right\} $, and $\left\{ \tau_{1},0,1\right\} $
in sequence to give $\frac{4\pi}{\sqrt{\eta_{2}^{2}-\eta_{1}^{2}}}\sinh^{-1}\left(\sqrt{\frac{\eta_{2}^{2}}{\eta_{1}^{2}}-1}\right)$.
But for the terms $n\geq1$, an infinite k-series would be required
so that this approach is not really an improvement on the spherical
coordinate equivalent (\ref{eq:theorem3}) and we will not discuss it further.

\section{Additional Coordinate Systems}

This one-range addition theorem may be applied to problems expressed
in other coordinates provided, of course, that the specific problem
results in a Slater orbital that contains a square root of coordinate
variables. 

\subsection{Ellipsoidal Coordinates}

Consider, for instance, the 90-year-old paper by Hirschfelder, Eyring,
and Rosen \cite{Hirschfelder} that calculates energy integrals related
to the $H_{3}$ molecule, whose ``nuclei \emph{a}, \emph{b}, and
\emph{c}, lie symmetrically on a straight line with a separation of
R times the first Bohr radius for atomic hydrogen $a_{0}$; and 1
and 2 are electrons.'' Their ellipsoidal coordinates are:

\begin{equation}
\begin{array}{ccc}
\lambda_{n} & = & \left(r_{in}+r_{jn}\right)/R'\\
\mu_{n} & = & \left(r_{in}-r_{jn}\right)/R'
\end{array}\label{eq:ellipsoidal}
\end{equation}
and $\phi_{n}$ is the azimuthal angle of n, ``where \emph{i} and
\emph{j} are the foci, $R'$ is the distance them, $r_{in}$ and
$r_{jn}$ are the distances from the arbitrary point \emph{n} to \emph{i}
and \emph{j}, respectively.'' In this coordinate system, one of the
several integrals they calculate has the form we desire to test,

\begin{equation}
T\left(a,bc\right)=2R^{3}\int_{1}^{\infty}d\lambda\int_{-1}^{1}d\mu\left(\frac{\lambda-\mu}{R}+\left(\lambda^{2}-\mu^{2}\right)\right)e^{-3R\lambda-R\mu}e^{-R\sqrt{\lambda^{2}+\mu^{2}-1}}\label{eq:T(a,bc)}
\end{equation}
except that the wave functions are hydrogenic rather than Slater orbitals.
This is no impediment if we have the following one-range addition
theorem at hand:

\textbf{Theorem 5.} For $k\leq1$,

\begin{equation}
\begin{array}{ccc}
e^{-x_{2}\sqrt{Bk^{2}+C}} & = & {\displaystyle \sum_{n=0}^{\infty}}\sqrt{{\displaystyle \frac{2}{\pi}}}{\displaystyle \frac{\left(-1\right)^{n}B^{n}k^{2n}}{n!}}2^{-n}x_{2}^{n+1/2}C^{\frac{1}{4}-\frac{n}{2}}K_{n-\frac{1}{2}}\left(\sqrt{C}x_{2}\right)\end{array}\label{eq:theorem5}
\end{equation}
 \textbf{Proof of Theorem 5.}

\begin{comment}
From the file ``Energy of H3 Molecule m7.nb''
\end{comment}
If we take the derivative of both sides of (\ref{eq:theorem1}) with
respect to $x_{2}$ and use the recursion relation \cite{functions.wolfram.com/03.04.17.0002.01}
the result immediately follows.$\square$ 

For the present problem, we insert (\ref{eq:theorem5}) in (\ref{eq:T(a,bc)})
and rewrite the half-integer Macdonald function as the finite series.
(\ref{eq:half-integer Macdonald function}) Then

\begin{align}
T\left(a,bc\right) & =2R^{3}\int_{1}^{\infty}d\lambda\int_{-1}^{1}d\mu\sum_{n=0}^{\infty}\sum_{J=0}^{\left\lfloor \left|n-\frac{1}{2}\right|-\frac{1}{2}\right\rfloor }\frac{2^{-J}R^{-J-\frac{1}{2}}\lambda^{-J-\frac{1}{2}}\left(\left|n-\frac{1}{2}\right|+J-\frac{1}{2}\right)!}{J!\left(\left|n-\frac{1}{2}\right|-J-\frac{1}{2}\right)!}\label{eq:T(a,bc)K_Series}\\
\times & \frac{\left(-\frac{1}{2}\right)^{n}R^{n+\frac{1}{2}}\left(\mu^{2}-1\right)^{n}e^{-4R\lambda-R\mu}}{n!}\left(\frac{\lambda^{\frac{3}{2}-n}-\mu\lambda^{\frac{1}{2}-n}}{R}-\mu^{2}\lambda^{\frac{1}{2}-n}+\lambda^{\frac{5}{2}-n}\right)\nonumber \\
= & \sum_{n=0}^{\infty}\int_{-1}^{1}\,d\mu\sum_{J=0}^{\left\lfloor \left|n-\frac{1}{2}\right|-\frac{1}{2}\right\rfloor }\frac{(-1)^{n}2^{J+n-5}R^{2n}\left(\mu^{2}-1\right)^{n}e^{-R\mu}\left(\left|n-\frac{1}{2}\right|+J-\frac{1}{2}\right)!}{J!n!\left(\left|n-\frac{1}{2}\right|-J-\frac{1}{2}\right)!}\nonumber \\
\times & (-16R\mu(R\mu+1)\Gamma(-J-n+1,4R)+4\Gamma(-J-n+2,4R)+\Gamma(-J-n+3,4R))\nonumber \\
= & \sum_{n=0}^{\infty}\sum_{J=0}^{\left\lfloor \left|n-\frac{1}{2}\right|-\frac{1}{2}\right\rfloor }\frac{i\sqrt{\pi}(-1)^{2n}e^{i\pi n}2^{j+2n-\frac{9}{2}}(-R)^{-n-\frac{3}{2}}R^{2n}\Gamma(n+1)\left(\left|n-\frac{1}{2}\right|+j-\frac{1}{2}\right)!\left(\right)}{j!n!\left(\left|n-\frac{1}{2}\right|-j-\frac{1}{2}\right)!}\nonumber \\
\times & \left(R\,I_{n+\frac{5}{2}}(R)\left(16R^{2}\Gamma(-j-n+1,4R)-4\Gamma(-j-n+2,4R)-\Gamma(-j-n+3,4R)\right)\right.\nonumber \\
- & \left.(2n+3)I_{n+\frac{3}{2}}(R)(4\Gamma(-j-n+2,4R)+\Gamma(-j-n+3,4R))\right)\:,\nonumber 
\end{align}
where the $I_{j}(R)$ are Modified Bessel functions and $\left\lfloor \left|n-\frac{1}{2}\right|-\frac{1}{2}\right\rfloor $
is the greatest integer less than $\left|n-\frac{1}{2}\right|-\frac{1}{2}$,
giving $0$ for $n=0$ and $n-1$ for higher values. Hirschfelder,
Eyring, and Rosen \cite{Hirschfelder} give the exact result as

\begin{equation}
\begin{array}{ccc}
T\left(a,bc\right) & = & \frac{1}{81}\left(e^{3R}\left(-16R^{2}+44R+\frac{116}{9R}-\frac{116}{3}\right)\text{Ei}(-8R)\right.\\
 & - & -e^{-3R}\left(16R^{2}+44R+\frac{116}{9R}+\frac{116}{3}\right)(\text{Ei}(-2R)+2\log(2))+\\
 & + & \left.\frac{1}{16}e^{-3R}\left(624R^{2}+2256R+\frac{131}{3R}+1670\right)-\frac{1}{16}e^{-5R}\left(160R+\frac{131}{3R}+34\right)\right)\:,
\end{array}\label{eq:exact}
\end{equation}
which for $R=.11$ gives $0.360071$. We find that the first six terms
of our one-range addition theorem give $0.356284+0.003537+0.00019+0.000036+0.000013+0.000005$
, and sum to $0.360061$. The seventh terms and beyond are dominated
by the $J=n-1$ term in the Macdonald function expansion, each giving
a contribution of $10^{-7}$. For $R=.011$ the convergence is much
more rapid but again plateaus, with the $J=n-1$ term in the Macdonald
function expansion each giving a contribution of $10^{-10}$. For
$R=1.1$ the convergence stalls in the fourth significant digit, with
the $J=n-1$ term in the Macdonald function expansion each giving
a contribution of $10^{-6}$. This behavior is very different from
the uniformly rapid convergence of the spherical coordinate example
of (\ref{eq:theorem3}), though whether it is a matter of the different
problem, or coordinate system, or choice of \emph{B} and \emph{C},
one could not say without attempting further examples in each system.
Given that for the arbitrary values $\left\{ C\to0.11,B\to0.13,x_{2}\to0.17,k\to0.23\right\} $
the first four terms in the one-range addition theorem itself (\ref{eq:theorem5}),
${0.945177-0.001665+0.000028-0.000000086}$ converge rapidly to the
value on the left-hand side, $0.943538$, one suspects that this behavior
has to do with the particular problem rather than use of this one-range
addition theorem within ellipsoidal coordinates more generally.

\subsection{Hylleraas coordinates}

Given the widespread use of Hylleraas coordinates \cite{Hyllerass1929}
for helium-like atoms and ions, 
\begin{equation}
\psi_{H}(r_{1},r_{2},r_{12})=\frac{1}{\sqrt{2}}\left(1-\hat{P}_{12}\right)e^{-\alpha r_{1}-\beta r_{2}-\gamma u}\sum_{l,m,n}c_{lmn}s^{l}t^{2m}u^{n},\label{eq:psi_H}
\end{equation}
where $\hat{P}_{12}$ is the permutation operator for the two identical
electrons and $s=r_{1}+r_{2}$, $t=r_{1}-r_{2}$, and $u=r_{12}\equiv\left|\mathbf{r}_{1}-\mathbf{r}_{2}\right|,$
some researchers might find use in a generalized version of the sequence
beginning with (\ref{eq:theorem1}) and (\ref{eq:theorem5}). Thus,
we offer the following infinite set of one-range addition theorems: 

\textbf{Theorem 6.} For $k\leq1$ and $j\geq0$,

\begin{equation}
\begin{array}{ccc}
\left(Bk^{2}+C\right)^{\frac{j-1}{2}}e^{-x_{2}\sqrt{Bk^{2}+C}} & = & {\displaystyle \sum_{n=0}^{\infty}}{\displaystyle \frac{1}{\sqrt{\pi}}\frac{\left(-1\right)^{n}B^{n}k^{2n}}{n!}}C^{\frac{j}{2}-n-\frac{1}{2}}G_{3,1}^{0,3}\left(\frac{4}{Cx_{2}^{2}}|\begin{array}{c}
\frac{1}{2},1,\frac{1}{2}(j-2n+1)\\
\frac{j+1}{2}
\end{array}\right)\:.\end{array}\label{eq:theorem6}
\end{equation}
 \textbf{Proof of Theorem 6.}

The proof structure follows that of (\ref{eq:theorem1}) except that
we use the more general Gaussian transform that is the second line
of (\ref{seventeen}) \cite{GR5}(p. 846 No. 7.386) and we require
the much more general inverse Gaussian transform \cite{PBM4}(p. 452
No. 3.25.2.7 with $k=l=1$ ), since our series expansion of $e^{-\left(Bk^{2}\right)\rho}$
gives $\rho$-powers of \emph{n} in addition to $\rho$-powers of
\emph{j} from (\ref{seventeen}), 

\begin{align}
\int_{0}^{\infty}\,d\rho\rho^{\mu}e^{-\frac{a^{2}}{\rho}-p\rho}H_{j}\left(\frac{a}{\sqrt{\rho}}\right) & =2^{j}p^{-\mu-1}G_{3,1}^{0,3}\left(\frac{1}{a^{2}p}|\begin{array}{c}
\frac{1}{2},1,-\mu\\
\frac{j+1}{2}
\end{array}\right)\:.\label{eq:p. 452 No. 3.25.2.7}\\
 & \left[Re\,p>0,\,\left|\arg a\right|<\pi/4\right]\nonumber 
\end{align}
In principle one could instead use \cite{PBM4}(p. 452 No. 3.25.2.5)
or \cite{PBM2} (p. 488 No. 2. 20.3. 22) that express this integral
in terms of gamma functions times $\,_{1}F_{2}$ hypergeometric functions
instead of the Meijer G-function, above, but in the present application
infinities sometimes ensue when the arguments of the gamma functions
become negative integers.$\square$ For the arbitrary values $\left\{ C\to0.11,\,B\to0.13,\,x_{2}\to0.17,\,k\to0.23\right\} $
and $j=0,\,1$ the results in the series (\ref{eq:theorem6}) indeed
precisely match those for (\ref{eq:theorem1}) and (\ref{eq:theorem5}). For
$j=2$, the first four terms in the series (\ref{eq:theorem6}) ${0.31348+0.00924-0.000161+0.000005}$
converge rapidly to the value on the left-hand side, $0.32257$.

\section*{Conclusion}

We have crafted an infinite set of one-range addition theorems for
Slater orbitals, and their derivatives, devoid of the infinite second
series that is typical of prior one-range addition theorems. Because
its derivation differs markedly from that for prior one-range addition
theorems, this approach may also be useful for other sorts of functions.
A test with the particular case of integrating over a product of Slater
orbitals (one having a center shifted from the origin) required an
additional infinite series to bring the problem to closure as a series
of analytic functions. Convergence is roughly one digit for every
other term in the series. When the charges on the Slater orbitals
were identical, this second series was unneeded. 

Unlike previous addition theorems, this set is applicable more than
one coordinate system. A particular example using ellipsoidal coordinates
gives a result in terms of a single infinite series, though the convergence
stalls out for this particular problem. Finally, we present an infinite
set of one-range addition theorems that can be used in Hylleraas coordinates. 

The one-range addition theorem for Slater orbitals can also be applied
to other Yukawa-like functions that may appear late in the reduction
of quantum amplitude integrals that include plane waves, allowing
the reduction of the final integral to an analytic series that converges
when the momentum variable $k\leq1$. The general result reduces to
a product of Macdonald and Whittaker functions for Cheshire's \cite{Cheshire}
integral that contains a product of a hydrogenic wave function, a
Slater orbital of the same charge, and a plane wave.


\begin{thebibliography}{2222222}
\bibitem[1]{Ley-Koo and Bunge}Ley-Koo, E.; Bunge, C.F. General evaluation
of atomic electron-repulsion integrals in orbital methods without
using a series representation for $r_{12}^{-1}$. \emph{Phys. Rev.
A} \textbf{1989}, \emph{40}, 1215.

\bibitem[2]{Sack}Sack, R.A. Generalization of Laplace's Expansion
to Arbitrary Powers and Functions of the Distance between Two Points.
\emph{J. Math. Phys.} \textbf{1964}, \emph{5}, 245.

\bibitem[3]{Porras and King}Porras, I.; King, F.W. Evaluation of
some integrals for the atomic three-electron problem using convergence
accelerators. \emph{Phys. Rev. A} \textbf{1994}, \emph{49}, 1637.

\bibitem[4]{Weniger_two-range}Weniger, E.J. Addition theorems as
three-dimensional Taylor expansions. II. B functions and other exponentially
decaying functions. \emph{ Int. J. Quant. Chem.} \textbf{2002}, \emph{90},
92\textendash 104.

\bibitem[5]{Guseinov2005}Guseinov I.I., One-range addition theorems
for combined Coulomb and Yukawa like central and noncentral interaction
potentials and their derivatives, \emph{J. Math. Chem.}\textbf{ 2006}
39 253\textendash 8.

\bibitem[6]{Fromm and Hill}Fromm, D.M.; Hill, R.N. Analytic evaluation
of three-electron integrals. \emph{Phys. Rev. A} \textbf{1987}, \emph{36},
1013.

\bibitem[7]{Remiddi}Remiddi, E. Analytic value of the atomic three-electron
correlation integral with Slater wave functions. \emph{Phys. Rev.
A} \textbf{1991}, \emph{44}, 5492.

\bibitem[8]{Harris PRA 55 1820}Harris, F.E. Analytic evaluation of
three-electron atomic integrals with Slater wave functions. \emph{Phys.
Rev. A} \textbf{1997}, \emph{55}, 1820.

\bibitem[9]{Kikuchi} Kikuchi, R. Gaussian Functions in Molecular
Integrals. \emph{J. Chem. Phys.} \textbf{1954}, \emph{22}, 148.

\bibitem[10]{Shavitt and Karplus}Shavitt, I.; Karplus, M. Multicenter
Integrals in Molecular Quantum Mechanics. \emph{ J. Chem. Phys.} \textbf{1962},
\emph{36}, 550.

\bibitem[11]{Stra89a} Straton, Jack C. Analytically reduced form
of multicenter integrals from Gaussian transforms. \emph{Phys. Rev.
A} \textbf{1989}, \emph{39}, 1676\textendash 1684; Errata: \emph{ Phys.
Rev. A} \textbf{1989} \emph{40}, 2819.

\bibitem[12]{stra23}Straton, Jack C. An integral representation for
quantum amplitudes. \emph{Phys. Scr. } \textbf{2023}, \emph{98}, 105406.

\bibitem[13]{stra24a}Straton, Jack C. Integral Representations over
Finite Limits for Quantum Amplitudes. Axioms \textbf{2024}, 13, 120.
%https://doi.org/10.3390/ axioms13020120.

\bibitem[14]{Chen} C. J. Chen, \emph{Introduction to Scanning Tunneling
Microscopy}, Oxford University Press, New York, Oxford, 1993, (Oxford
Series in Optical and Imaging Science 4, Eds. M. Lapp, J.-I. Nishizawa,
B. B. Snavely, H. Stark, A. C. Tam, T. Wilson, ISBN 0-19-507150-6)
p. 151, eq. (6.8).

\bibitem[15]{Yukawa} Yukawa, H. On the Interaction of Elementary
Particles. I. \emph{Proc. Phys. Math. Soc. Jpn.} \textbf{1935}, \emph{17},
48.

\bibitem[16]{NayekGhoshal} Nayek, S.; Ghoshal, A. Dynamics of positronium
formation in positron-hydrogen collisions embedded in weakly coupled
plasmas. \emph{Phys. Plasmas} \textbf{2012}, \emph{19}, 113501.



\bibitem[17]{Harris}Harris, G.M. Attractive Two-Body Interactions
in Partially Ionized Plasmas. \emph{Phys. Rev.} \textbf{1962}, \emph{125},
1131.


\bibitem[18]{EckerWeizel} Ecker, G.; Weizel, W. Zustandssumme und
effective Ionisierungsspannung eines Atoms im Innern des Plasmas.
\emph{Ann. Phys.} \textbf{1956}, \emph{17}, 126.

\bibitem[19]{Smirnov2003}B. M. Smirnov, \emph{Physics of Atoms and
Ions} (Springer-Verlag, New York, 2003), p. 190.

\bibitem[20]{GaravelliOliveira} Garavelli, S.L.; Oliveira, F.A. Analytical
solution for a Yukawa-type potential. \emph{Phys. Rev. Lett.} \textbf{1991},
\emph{66}, 1310.

\bibitem[21]{CaccavanoLeung} Caccavano, A.; Leung, P.T. Atomic spectroscopy
and the photon mass: Effects on the 21 cm radiation. \emph{Phys. Lett.
A} \textbf{2013}, \emph{377}, 2777.

\bibitem[22]{Stilgoe}Stilgoe, A. B., Nieminen, T. A. and Rubinsztein-Dunlop,
H. Computational toolbox for scattering of focused light from flattened
or elongated particles using spheroidal wavefunctions. Journal of
Quantitative Spectroscopy and Radiative Transfer, \textbf{2025} \emph{331}
109267. %doi: 10.1016/j.jqsrt.2024.109267 

\bibitem[23]{Joachain} C. J. Joachain, \emph{Quantum Collision Theory}; North-Holland, NY, 1983.

%\bibitem[23]{Goldberger and K. M. Watson} M. L. Goldberger and K. M. Watson, \textbackslash{}textit\{Collision Theory\} (John Wiley, NY, 1964) .

\bibitem[24]{Casey}Yazejian, C. A.  and Straton, Jack C.  Polarization
in the production of the antihydrogen ion, \emph{Euro. Phy. J. D}
\textbf{2020}, \emph{74}, 156. % DOI: 10.1140/epjd/e2020-100548-7.

\bibitem[25]{Magnus Oberhettinger and Soni}W. Magnus, F. Oberhettinger,
and R. P. Soni, F\emph{ormulas and Theorems for the Special Functions
of Mathematical Physics} (Springer - Verlag : NY, 1966) . p. 107 (the
1954 edition, p 21 has the i k version).

\bibitem[26]{GR5} Gradshteyn, I.S.; Ryzhik, I.M. \textit{Table of
Integrals, Series, and Products}, 5th ed.; Academic: \mbox{New York,
NY, USA,} 1994.

\bibitem[27]{Guseinov2002} Guseinov I.I., New Complete Orthonormal
Sets of Exponential-Type Orbitals and Their Application to Translation
of Slater Orbitals, \emph{Int. J. Quantum Chem.} \textbf{2002}, \emph{90},
114.

\bibitem[28]{Cheshire} {Cheshire, I.M. Positronium formation by
fast positrons in atomic hydrogen}. \emph{Proc. Phys. Soc.} \textbf{{1964}},
{83}, {227\textendash 237.}

\bibitem[29]{functions.wolfram.com/06.34.26.0002.01} Available online:
https://functions.wolfram.com/06.34.26.0002.01 (accessed on 9 January
2025).

\bibitem[30]{functions.wolfram.com/06.06.26.0002.01} Available online:
https://functions.wolfram.com/06.06.26.0002.01 (accessed on 9 January
2025).

\bibitem[31]{functions.wolfram.com/06.25.26.0001.01} Available online:
https://functions.wolfram.com/06.25.26.0001.01 (accessed on 9 January
2025).

\bibitem[32]{mathworld.wolfram.com/LegendrePolynomial.html}Weisstein,
Eric W. \textquotedbl Legendre Polynomial.\textquotedbl{} From MathWorld\textendash A
Wolfram Resource. Available online: https://mathworld.wolfram.com/LegendrePolynomial.html
(accessed on 9 January 2025).

\bibitem[33]{PBM5}Prudnikov, A.P.; Brychkov, Y.A.; Marichev, O.I.
\textit{Integrals and Series}; Gordon and Breach: \mbox{New York,
NY, USA}, 1986; Volume 5.

\bibitem[34]{PBM3}Prudnikov, A.P.; Brychkov, Y.A.; Marichev, O.I.
\textit{Integrals and Series}; Gordon and Breach: \mbox{New York,
NY, USA}, 1986; Volume 3.

\bibitem[35]{Hirschfelder}Hirschfelder, J., Eyring, H., and Rosen,
N. I. Calculation of Energy of H3 Molecule, J. Chem. Phys. \textbf{1936}
\emph{4}, 121\textendash 130. %https://doi.org/10.1063/1.1749798

\bibitem[36]{functions.wolfram.com/03.04.17.0002.01}Available online:
http://functions.wolfram.com/03.04.17.0002.01 (accessed on 9 January
2025).

\bibitem[37]{Hyllerass1929}Hyllerass, E. A. Z. Physik \textbf{1929},
\emph{54}, 347.

\bibitem[38]{PBM4}Prudnikov, A.P.; Brychkov, Y.A.; Marichev, O.I.
\textit{Integrals and Series}; Gordon and Breach: \mbox{New York,
NY, USA}, 1992; Volume 4.

\bibitem[39]{PBM2}Prudnikov, A.P.; Brychkov, Y.A.; Marichev, O.I.
\textit{Integrals and Series}; Gordon and Breach: \mbox{New York,
NY, USA}, 1986; Volume 2.

%\bibitem[1888888]{Slater}Slater, L.J. Expansions of Generalized Whittaker
%Functions. \emph{Math. Proc. Camb. Philos. Soc.} \textbf{1954}, {\em
%50}, 628--631.
%\bibitem[333333]{Stra90a}Jack C. Straton, Reduced form for the general-state
%multicenter integral from an integro-differential transform. \emph{Phys.
%Rev. A} \textbf{1990}, \textit{41}, 71--77. 

\end{thebibliography}
\end{document}